\documentclass[11pt]{article}

\usepackage[T1]{fontenc}
\usepackage[utf8]{inputenc}
\usepackage[english]{babel}
\usepackage{amsmath,amssymb,amsfonts}
\usepackage{indentfirst}
\usepackage[margin=1in]{geometry}
\usepackage[hidelinks]{hyperref}
\hypersetup{
  pdftitle={Sharp Bernstein--Nikolskii inequalities for polynomials and entire functions of exponential type},
  pdfauthor={D. V. Gorbachev}
}

\renewcommand{\le}{\leqslant}
\renewcommand{\ge}{\geqslant}
\title{Sharp Bernstein--Nikolskii inequalities for polynomials and entire functions of exponential type\thanks{This work was supported by the Russian Foundation for Basic Research (RFBR), project No.~20-11-50107.}}
\author{D.~V.~Gorbachev}
\date{}

\begin{document}
\maketitle

\begin{abstract}
The classical Bernstein--Nikolskii inequalities of the form $\|Df\|_{q}\le
\mathcal{C}_{pq}\|f\|_{p}$, $f\in Y$, estimate the $pq$-norms of differential
operators $D$ on classes $Y$ of polynomials and entire functions of exponential
type. These inequalities play an important role in harmonic analysis and
approximation theory and also have applications in number theory and metric
geometry. Both order-sharp inequalities and inequalities with sharp constants
are studied. The latter are especially interesting because extremal functions
reflect the geometry of the underlying manifold, which can be useful in
geometric problems.

Historically, Bernstein inequalities correspond to the case $p=q$, while
Nikolskii inequalities concern the identity operator for $p<q$. The first
estimate for the derivative of a trigonometric polynomial for $p=\infty$ was
obtained by S.N.~Bernstein (1912), although A.A.~Markov had earlier given its
algebraic analogue (1889). Bernstein's inequality was refined by E.~Landau and
M.~Riesz, and A.~Zygmund proved it for all $p\ge1$ (1933). For $p<1$, an
order-sharp Bernstein inequality was obtained independently by V.I.~Ivanov
(1975) and by E.A.~Storozhenko, V.G.~Krotov, and P.~Oswald (1975), while the
sharp inequality was proved by V.V.~Arestov (1981). For entire functions of
exponential type, the sharp Bernstein inequality was proved by N.I.~Akhiezer
and B.Ya.~Levin for $p\ge1$ (1957), and by Q.I.~Rahman and G.~Schmeisser for
$p<1$ (1990).

The first one-dimensional Nikolskii inequalities for $q=\infty$ were obtained
by D.~Jackson (1933) for trigonometric polynomials and by J.~Korevaar (1949)
for entire functions of exponential type. The general $d$-dimensional result
for $q\le\infty$ was established by S.M.~Nikolskii (1951). Estimates for
Nikolskii constants were later refined by I.I.~Ibragimov (1959), D.~Amir and
Z.~Ziegler (1976), R.J.~Nessel and G.~Wilmes (1978), and many others.
Order-sharp Bernstein--Nikolskii inequalities on different intervals were
studied by N.K.~Bari (1954). Versions for general multiplier differential
operators and weighted settings can be found in works of P.I.~Lizorkin (1965),
A.I.~Kamzolov (1984), A.G.~Babenko (1992), A.I.~Kozko (1998), K.V.~Runovskii
and H.-J.~Schmeisser (2001), F.~Dai and Y.~Xu (2013), V.V.~Arestov and
P.Yu.~Glazyrina (2014), and other authors.

For a long time, Bernstein--Nikolskii inequalities for polynomials and for
entire functions of exponential type developed as parallel theories. In 2015,
E.~Levin and D.~Lubinsky proved that for every $p>0$ the Nikolskii constant for
entire functions is the limit of the corresponding trigonometric constants.
For Bernstein--Nikolskii constants, this was proved by M.I.~Ganzburg and
S.Yu.~Tikhonov (2017) and refined by the author and I.A.~Martyanov (2018,
2019). Multidimensional Levin--Lubinsky type results were obtained by the
author jointly with F.~Dai and S.Yu.~Tikhonov for the sphere (2020), and by
M.I.~Ganzburg for the torus (2019) and the cube (2021).

At present, sharp Nikolskii constants are known in general only for
$(p,q)=(2,\infty)$. The case $p=1$ is particularly intriguing. Progress in
this problem was made by Ya.L.~Geronimus (1938), S.B.~Stechkin (1961),
L.V.~Taikov (1965), L.~H\"ormander and B.~Bernhardsson (1993), N.N.~Andreev,
S.V.~Konyagin, and A.Yu.~Popov (1996), the author (2005), the author and
I.A.~Martyanov (2018), and I.E.~Simonov and P.Yu.~Glazyrina (2015).
E.~Carneiro, M.B.~Milinovich, and K.~Soundararajan (2019) found applications to
the Riemann zeta function. V.V.~Arestov, M.V.~Deikalova, and their coauthors
(2016, 2018) characterized extremal polynomials for general weighted Nikolskii
constants using duality; earlier work in this direction goes back to
S.N.~Bernstein, L.V.~Taikov, and others.

A newer direction is the study of sharp Nikolskii inequalities on constrained
classes of functions. This reveals connections with extremal problems of
harmonic analysis, including the problems of Tur\'an and Delsarte and the
uncertainty principle of J.~Bourgain, L.~Clozel, and J.-P.~Kahane (2010). For
example, the author and coauthors (2020) showed that the sharp Nikolskii
constant for nonnegative spherical polynomials yields the
Delsarte--Goethals--Seidel bound for spherical designs. Related problems for
functions lead to well-known bounds for sphere-packing density, while
order-sharp results are closely connected with Fourier inequalities.

These results are presented within a common framework for Bernstein--Nikolskii
inequalities. We also discuss applications in approximation theory, number
theory, and metric geometry, and formulate several open problems.
\end{abstract}


\section{Notation}

\begin{itemize}

\item $\mathbb{Z}_{+}=\mathbb{Z}_{\ge 0}$, $\mathbb{R}_{+}=\mathbb{R}_{\ge 0}$;

\item $L^{p}(M)$~--- the space of measurable functions $f\colon M\to \mathbb{C}$ with
finite (quasi-)norm $\|f\|_{p}=\|f\|_{p;M}=(\int_{M}|f|^{p}\,dx)^{1/p}$ for
$0<p<\infty$ and $\|f\|_{\infty}=\operatorname*{ess\,sup}_{M}|f|$ for $p=\infty$;

\item $C(M)$~--- the space of continuous functions,
$\|f\|_{C}=\|f\|_{\infty}=\sup_{M}|f|$;

\item $L_{w}^{p}(M)$~--- the weighted space,
$\|f\|_{p,w}=(\int_{M}|f|^{p}w\,dx)^{1/p}$ for $p<\infty$,
$w\colon M\to \mathbb{R}_{+}$~--- a weight;

\item $p'=\frac{p}{p-1}$ for $1<p<\infty$, $1'=\infty$~--- the conjugate exponent;

\item $\mathbb{T}=(-\pi,\pi]$~--- the one-dimensional torus;

\item $\mathbb{R}^{d}$~--- the $d$-dimensional Euclidean space,
$xy=x_{1}y_{1}+\ldots+x_{d}y_{d}$~--- the scalar product of $x,y\in
\mathbb{R}^{d}$, $|x|=(xx)^{1/2}$;

\item $B^{d}\subset \mathbb{R}^{d}$~--- the unit Euclidean ball centered at the
origin;

\item $\mathbb{S}^{d}\subset \mathbb{R}^{d+1}$~--- the unit Euclidean sphere;

\item $\chi_{A}$~--- the characteristic function of a set $A$;

\item $|A|$~--- the number of elements of a finite set $A$ or the measure of a
measurable set $A$;

\item $\operatorname{supp}f$~--- the support of a function $f$;

\item $\widehat{f}(y)=(2\pi)^{-d}\int_{\mathbb{R}^{d}}f(x)e^{-ixy}\,dx$~---
the Fourier transform of $f$;

\item $(f*g)(x)=\int_{\mathbb{R}^{d}}f(x-y)g(y)\,dy$~--- the convolution of
$f$ and $g$;

\item $C$~--- positive constants that may vary from line to line and depend
only on inessential parameters;

\item $A(\theta)\asymp B(\theta)$ if $C_{1}A(\theta)\le B(\theta)\le
C_{2}A(\theta)$, where $C_{1}\le C_{2}$ are independent of $\theta$;

\item $(a)_{+}=\max\,\{a,0\}$;

\item $[a]$~--- the integer part, while $\lceil a\rceil$ and $\lfloor a\rfloor$~---
rounding up and down, respectively;

\item $Y\subset C(M)$~--- usually a subspace;

\item $K$~--- the reproducing kernel of the subspace $Y$;

\item $\rho$~--- the metric of a homogeneous space $M=G/H$,
$G$~--- a transitive group of motions of $M$, $H$~--- the stabilizer of a point
$o\in M$, $f(x)=f(\rho(x,o))$~--- zonal functions;

\item $f(x)=f(|x|)$~--- radial functions on $\mathbb{R}^{d}$;

\item $D$~--- a differential operator, $I$~--- the identity operator;

\item $\partial=\partial_{x}$~--- differentiation with respect to $x$;

\item $\Delta_{M}$~--- the Laplace--Beltrami operator on the homogeneous space
$M$, $\Delta_{0}=\Delta_{\mathbb{S}^{d}}$;

\item $\mathcal{C}_{pq}^{D}(Y)$~--- a Bernstein--Nikolskii constant,
$\mathcal{C}_{pq,w}^{D}(Y)$~--- a weighted constant,
$\mathcal{C}_{p}^{D}(Y)=\mathcal{C}_{p\infty}^{D}(Y)$ when $q=\infty$;

\item $\mathcal{A}_{p}^{D}(Y;o)$~--- the shifted Nikolskii problem;

\item $\mathcal{T}_{n}$~--- the set of complex-valued $2\pi$-periodic
trigonometric polynomials of order at most $n$;

\item $\mathcal{P}_{n}$~--- the set of complex-valued algebraic polynomials of
degree at most $n$;

\item $\mathcal{E}_{\sigma}$~--- the set of entire functions of exponential type
at most $\sigma$, and $\mathcal{E}_{\sigma}^{p}$~--- the subset of functions with
finite $L^{p}$-norm;

\item $\mathbf{T}_{n}(x)$, $\mathbf{U}_{n}(x)$~--- the Chebyshev polynomials of
the first and second kind, respectively;

\item $(1-x)^{\alpha}(1+x)^{\beta}$, $P_{n}^{(\alpha,\beta)}(x)$~--- the Jacobi
weight and Jacobi polynomials, respectively; for $\alpha=\beta$ they are the
Gegenbauer weight and polynomials;

\item $j_{\alpha}(x)=\frac{2^{\alpha}\Gamma(\alpha+1)J_{\alpha}(x)}{x^{\alpha}}$~---
the normalized Bessel function;

\item $0<q_{\alpha 1}<q_{\alpha 2}<\ldots$~--- the zeros of the Bessel function
$J_{\alpha}(x)$;

\item $\operatorname{Si}x=\int_{0}^{x}\frac{\sin t}{t}\,dt$~--- the sine
integral.

\end{itemize}

\section{Introduction}

Bernstein--Nikolskii inequalities are studied on subsets of the Lebesgue--Riesz
spaces $L^{p}$ of real- or complex-valued functions defined on a manifold~$M$.
For the natural measure $dx$ on $M$, we write
$\|f\|_{p}=(\int_{M}|f|^{p}\,dx)^{1/p}$ for $0<p<\infty$ (a quasi-norm when
$p<1$), and $\|f\|_{\infty}=\operatorname*{ess\,sup}_{M}|f|$ for $p=\infty$.
For continuous functions in $C(M)$ we use the uniform norm
$\|f\|_{\infty}=\sup_{M}|f|$. As usual, for positive functions $A(\theta)$ and
$B(\theta)$, the notation $A(\theta)\asymp B(\theta)$ means that
$C_{1}A(\theta)\le B(\theta)\le C_{2}A(\theta)$ with constants
$0<C_{1}\le C_{2}$ independent of the parameter $\theta$.

Classical Bernstein--Nikolskii inequalities have the form
\[
\|Df\|_{q}\le C\|f\|_{p},\quad f\in Y,
\]
where $D$ is a differential operator (including the identity $I$), and $Y$ is
usually a subspace of $C(M)$ (in our setting, a space of polynomials or entire
functions of exponential type), with $DY\subset Y$. The sharp constant in this
inequality is defined as the solution of the Bernstein--Nikolskii problem
\[
\mathcal{C}_{pq;M}^{D}(Y)=\sup_{f\in Y\cap L^{p}(M)}
\frac{\|Df\|_{q}}{\|f\|_{p}}
\]
(with $0/0$ understood as $0$). For brevity, we omit $M$ and $D=I$ from the
notation when no confusion can arise. This quantity is the $pq$-operator norm
of $D$ on the subspace~$Y$. Note that $Y$ itself may be larger than
$Y\cap L^{p}(M)$, but only the latter set enters the supremum.

It is natural to study $\mathcal{C}_{pq}^{D}(Y)$ as a function of all its
parameters, and especially of the size of $Y$, for example of $\dim Y$. More
general settings allow nonlinear classes $Y$, fractional operators $D$, and
other modifications.

The cases $D=I$, $p<q$ and $D\ne I$, $p=q$ correspond to the classical
Nikolskii and Bernstein inequalities, respectively. Bernstein-type inequalities
for algebraic polynomials are also called Markov inequalities. A historical
overview explaining this terminology is given in Section~\ref{sec-hist}.

Bernstein--Markov--Nikolskii inequalities and extremal problems are widely used
in approximation theory to prove direct and inverse theorems for best
approximation by subspaces~$Y$. Similar extremal problems arise in analytic
number theory, coding theory, and metric geometry, especially when the classes
$Y$ are varied. Some of these applications are discussed in
Section~\ref{sec-appl}.

A detailed account of the classical Markov--Bernstein--Nikolskii inequalities,
mainly in one dimension, is given in the well-known monograph by
G.V.~Milovanovi\'c, D.S.~Mitrinovi\'c, and Th.M.~Rassias~\cite{MilMitRas94}.
Here we therefore focus on recent progress in multidimensional sharp
Bernstein--Nikolskii inequalities in the spaces $L^{p}(M)$ on homogeneous
manifolds $M$, including the unit Euclidean sphere $\mathbb{S}^{d}$ and the
Euclidean space~$\mathbb{R}^{d}$. Many of these results rely on solutions of
one-dimensional problems in weighted $L^{p}$ spaces, which enriches the
classical theory.

Computing $\mathcal{C}_{pq}^{D}(Y)$ in full generality is difficult. At
present, sharp results are essentially limited to the pairs
$(p,q)=(2,\infty)$ and $(p,p)$. We discuss them below. Many other results are
therefore only order estimates. In particular, let $D$ be a homogeneous
operator of order~$r$, and let $Y_{n}$ depend on a parameter $n$, for example
so that $\dim Y_{n}\asymp n^{d}$ as $n\to\infty$. A typical estimate for the
Bernstein--Nikolskii constant then has the form
\begin{equation}\label{C-asymp}
\mathcal{C}_{pq}^{D}(Y_{n})\asymp
n^{r+d(1/p-1/q)_{+}},\quad n\to \infty.
\end{equation}
Among the main results discussed in this survey are Levin--Lubinsky type
theorems, where this equivalence can be replaced by an asymptotic equality
\begin{equation}\label{ll-cL}
\mathcal{C}_{pq}^{D}(Y_{n})=\mathcal{L}_{pq}^{D}
 n^{r+d(1/p-1/q)_{+}}(1+o(1)),
\end{equation}
with a limiting Bernstein--Nikolskii constant $\mathcal{L}_{pq}^{D}$. The
theory of such asymptotic formulas is still developing, and at present they
are known only for $q=\infty$, which will be the main case in this survey.

The survey is organized as follows. Section~\ref{sec-hist} gives the historical
background, with emphasis on the results most relevant here.
Section~\ref{sec-common} presents, in a compact form, standard approaches to
upper and lower estimates for Bernstein--Nikolskii constants, together with
several related results. In Section~\ref{sec-ll} we prove explicit bounds for
Nikolskii constants associated with the multidimensional sphere
$\mathbb{S}^{d}$; these bounds imply a Levin--Lubinsky type result with an
explicit remainder estimate. We chose this result because it is fairly general
and its proof is comparatively simple. Section~\ref{sec-appl} discusses several
applications of Bernstein--Nikolskii constants and extremal problems, including
approximation theory, number theory, and metric geometry. Finally, we list some
open problems concerning sharp Bernstein--Nikolskii inequalities.

\section{Historical overview of the main results}\label{sec-hist}

We shall follow the historical development of the subject to some extent, while
selecting from the large body of results those that are conceptually important
for us. For convenience, we divide the discussion into periods according to
the main results of each period.

\subsection{The era of Chebyshev and the Markov brothers}
The development of classical inequalities in approximation theory for
polynomials and functions was largely stimulated by P.L.~Chebyshev's famous
problem of finding a polynomial of least deviation from zero. Our account of
the classical results follows the well-known book~\cite{MilMitRas94}. In what
follows, $\mathcal{P}_{n}$ denotes the set of algebraic polynomials
$P(x)=\sum_{k=0}^{n}c_{k}x^{k}$ with complex coefficients and degree at most
$n\in \mathbb{Z}_{+}$. Here and below, the subset of real polynomials is marked
by the superscript ``real''. Then, for $M=[-1,1]$ and every $n\ge 1$,
\begin{equation}\label{min-Pn}
\min_{P\in
\mathcal{P}_{n-1}^{\text{real}}}\|x^{n}+P(x)\|_{\infty}
=2^{1-n}\|\mathbf{T}_{n}(x)\|_{\infty},
\end{equation}
where $\mathbf{T}_{n}(x)=\cos n\arccos x$ is the unique extremal polynomial,
later called the Chebyshev polynomial of the first kind. It is believed that
this problem arose in Chebyshev's work on Watt's linkage. The method of its
solution later led to the famous alternation theorem, a criterion for best
approximation by a subspace in the uniform norm. By compactness, we write a
minimum in~\eqref{min-Pn} from the outset.

Chebyshev polynomials and their generalizations, including some rather
unexpected ones, appear in many problems (see, for example,~\cite{AndYud01}).
Here we single out a result of the Markov brothers, which states that the
Chebyshev polynomial is extremal in the problem of estimating derivatives of
algebraic polynomials. Namely, for every real polynomial
$P_{n}\in \mathcal{P}_{n}$ and $r\in \mathbb{N}$,
\begin{equation}\label{mar-ineq}
\|P_{n}^{(r)}\|_{\infty}\le \mathbf{T}_{n}^{(r)}(1)\|P_{n}\|_{\infty},
\end{equation}
or, writing $\partial=\partial_{x}$ for differentiation with respect to $x$,
\[
\mathcal{C}_{\infty\infty}^{\partial^{r}}(\mathcal{P}_{n})
=\mathbf{T}_{n}^{(r)}(1).
\]
For $r=1$, inequality~\eqref{mar-ineq} was proved by A.A.~Markov, and the case
$r\ge1$ was proved by his brother V.A.~Markov~\cite{MilMitRas94}. It is
interesting that this inequality grew out of D.I.~Mendeleev's sharp estimate
$\|P_{2}'\|_{\infty}\le4\|P_{2}\|_{\infty}$, which he used in his well-known
work on the concentration of alcohol in water (see~\cite{MilMitRas94}). Later,
extremal problems for derivatives of polynomials on sets
$M\subset\mathbb{R}^{d}$ came to be called Markov-type inequalities or
inequalities of the Markov brothers. V.A.~Markov also proved the following
coefficient version:
\begin{equation}\label{Pnk-0}
|P_{n}^{(r)}(0)|\le |\mathbf{T}_{n-\theta}^{(r)}(0)|\,\|P\|_{\infty},
\end{equation}
where $\theta=0$ or $1$ according as $n-r$ is even or odd. For $r=n$, this is
known as Chebyshev's inequality.

A.N.~Korkin and E.I.~Zolotarev solved Chebyshev's problem in the integral
metric of $L^{1}([-1,1])$~\cite{MilMitRas94}:
\begin{equation}\label{kz-res}
\min_{P\in
\mathcal{P}_{n-1}^{\text{real}}}\|x^{n}+P(x)\|_{1}
=2^{-n}\|\mathbf{U}_{n}(x)\|_{1},
\end{equation}
where
$\mathbf{U}_{n}(x)=\frac{\sin((n+1)\arccos x)}{\sqrt{1-x^{2}}}$ is the
Chebyshev polynomial of the second kind, which is also unique in this problem.

Chebyshev's problem is relatively simple in $L^{2}([-1,1])$, where the unique
solution is given by the orthogonal Legendre polynomials. For other $p\ge1$,
extremal polynomials can be characterized by orthogonality relations, but no
constructive description of these polynomials is known to us.

Chebyshev problems have also been studied extensively for polynomials with
several fixed coefficients (see~\cite{MilMitRas94}). We shall pay more
attention to weighted analogues of Chebyshev's problem, because they are
closely related to multidimensional sharp Nikolskii inequalities. Given a
nonnegative weight $\varpi$ on $[-1,1]$, one has to find
\[
\min_{P\in \mathcal{P}_{n-1}^{\text{real}}}
\|(x^{n}+P(x))\varpi(x)\|_{p}.
\]
In the $L^{2}$ metric the answer is again given by orthogonal polynomials, now
with weight $\varpi^{2}$. The answer is also known for $p=1,\infty$ for weights
of the form
\[
\varpi(x)=\frac{(1-x)^{a}(1+x)^{b}}{\sqrt{S(x)}},\quad
a,b\in \{0,1/2\},
\]
where $S$ is a polynomial positive on $[-1,1]$ (see
\cite{MilMitRas94,Leb04}). In all these cases the extremizers are variants of
the Chebyshev--Markov--Bernstein--Szeg\H{o} polynomials.

Of particular interest to us is the Jacobi weight
$\varpi(x)=(1-x)^{\alpha}(1+x)^{\beta}$~\cite{Seg62}. For general $\alpha$ and
$\beta$, however, Chebyshev's problem is not solved when $p\ne2$, although
order estimates are known~\cite{MilMitRas94}. Below we shall see how this case
is related to the Nikolskii inequality for spherical polynomials.

There are also interesting studies of the case where $\varpi\equiv0$ on
intervals. Then Chebyshev's problem becomes close to the Remez problem (see,
for example,~\cite{Bog99,TikYud20} and Subsection~\ref{subsec-appl-remez}). In
the $L^{\infty}$ Remez problem the extremal polynomials are again modifications
of the Chebyshev polynomial $\mathbf{T}_{n}$~\cite{TikYud20}.

After the trigonometric substitution $x=\cos t$, many algebraic problems reduce
to problems for even trigonometric polynomials in $L^{p}(\mathbb{T})$ on the
torus $\mathbb{T}=(-\pi,\pi]$. We denote by $\mathcal{T}_{n}$ the set of all
complex-valued trigonometric polynomials
$T(t)=\sum_{k=-n}^{n}c_{k}e^{ikt}$ of order at most $n$. Even and odd
polynomials will be marked by the superscripts ``even'' and ``odd''. Then the
algebraic Chebyshev problem is equivalent to the following trigonometric
problem in $L^{\infty}(\mathbb{T})$:
\[
\min_{T\in \mathcal{T}_{n-1}^{\text{real},\text{even}}}
\|\!\cos nt+T(t)\|_{\infty}=\|\!\cos nt\|_{\infty}.
\]
The trigonometric problem can be solved by a method different from alternation,
namely by averaging over the discrete subgroup of $\mathbb{T}$ consisting of
the points $t_{k}=\frac{2\pi k}{n}$, $k=0,1,\ldots,n-1$. With this method the
evenness condition can be dropped (see~\cite{AndYud01}). The same idea extends
to arbitrary translation-invariant norms. In particular, S.N.~Bernstein
proved~\cite{Tim60} that for every $p\in[1,\infty]$ in
$L^{p}(\mathbb{T})$,
\begin{equation}\label{T-cos}
\min_{T\in \mathcal{T}_{n-1}}
\|\!\cos n(t-\tau)+T(t)\|_{p}=\|\!\cos nt\|_{p}.
\end{equation}
Indeed, for every $T\in\mathcal{T}_{n-1}$ and $t\in\mathbb{T}$,
$\frac{1}{n}\sum_{k=0}^{n-1}T(t-t_{k})=c_{0}$. Hence, by the triangle
inequality, translation invariance of the norm, and the periodicity relation
$\cos n(t+t_{k}-\tau)=\cos n(t-\tau)$,
\[
\|\!\cos n(t-\tau)+c_{0}\|_{p}\le
\frac{1}{n}\sum_{k=0}^{n-1}
\|\!\cos n(t-\tau)+T(t-t_{k})\|_{p}
=\|\!\cos n(t-\tau)+T(t)\|_{p}.
\]
After the shift $t\mapsto t+\frac{\pi}{n}$ we obtain
\[
\|\!\cos n(t-\tau)-c_{0}\|_{p}\le
\|\!\cos n(t-\tau)+T(t)\|_{p}.
\]
A second application of the triangle inequality gives
$\|\!\cos n(t-\tau)\|_{p}\le
\|\!\cos n(t-\tau)+T(t)\|_{p}$, with equality for the zero polynomial.
This approach was extended to the multidimensional setting by N.N.~Andreev and
V.A.~Yudin~\cite{AndYud01}.

It follows from~\eqref{T-cos} that
\begin{equation}\label{T-sin}
\min_{T\in \mathcal{T}_{n}^{\text{real},\text{odd}}}
\|\!\sin((n+1)t)+T(t)\|_{1}=\|\!\sin((n+1)t)\|_{1},
\end{equation}
and the inverse trigonometric substitution yields the Korkin--Zolotarev
result~\eqref{kz-res}. So far we have formulated the problems mainly for real
polynomials, but, as we shall see, many results remain valid for complex-valued
polynomials as well.

\subsection{The era of Bernstein and Jackson}
Equality~\eqref{T-sin} is equivalent to the algebraic $L^{1}$ Chebyshev problem,
and conversely. Under the trigonometric substitution, the weight
$|\!\sin t|$ appears. For the same reason, the algebraic Markov problem for
derivatives is not equivalent to the corresponding problem for derivatives of
trigonometric polynomials. S.N.~Bernstein first proved such an inequality for
even trigonometric polynomials $T_{n}\in\mathcal{T}_{n}$ (see
\cite{MilMitRas94}): for $M=\mathbb{T}$,
\[
\|T_{n}^{(r)}\|_{\infty}\le n^{r}\|T_{n}\|_{\infty}.
\]
Later, through the work of E.~Landau ($p=\infty$, arbitrary real polynomials,
not necessarily even), M.~Riesz ($p=\infty$, complex polynomials), A.~Zygmund
($p\ge1$), and V.V.~Arestov ($0\le p<1$), this inequality was established in
its full generality~\cite{MilMitRas94}: for every $p\in[0,\infty]$,
\begin{equation}\label{ber-ineq}
\|T_{n}^{(r)}\|_{p}\le n^{r}\|T_{n}\|_{p},
\end{equation}
that is,
$\mathcal{C}_{pp}^{\partial^{r}}(\mathcal{T}_{n})=n^{r}$ (for the definition of
the $0$-quasi-norm, see~\cite{Are81}). The polynomials
$C\cos n(t-\tau)$ are extremal.

A.~Zygmund used the M.~Riesz interpolation formula
\begin{equation}\label{Tn-sum}
T_{n}'(t)=\sum_{j=1}^{2n}a_{j}T_{n}(t+t_{j}),\quad
\sum_{j=1}^{2n}|a_{j}|=n,\quad
t_{j}=\frac{(2j-1)\pi}{2n},\quad
a_{j}=\frac{(-1)^{j+1}}{4n\sin^{2}\frac{t_{j}}{2}},
\end{equation}
and H\"older's inequality, which is available for $p\ge1$. For $p<1$ this
argument cannot be used, and the problem is substantially harder. It was solved
with a nonsharp constant by V.I.~Ivanov~\cite{Iva75} and independently by
E.A.~Storozhenko, V.G.~Krotov, and P.~Oswald~\cite{StoKroOsv75}; see also
\cite{AttNev80}. The final solution of the classical Bernstein problem was
obtained by V.V.~Arestov~\cite{Are81} using methods of complex analysis. A
relatively simple version of the proof can be found in~\cite{Erd20}. In our
view, many open problems for real polynomials require going into the complex
domain. In this survey, however, we try to stay mainly within real-variable
methods.

Several proofs of the sharp Bernstein inequality~\eqref{ber-ineq} can be found
in~\cite{QueZar19}. Another proof, for $p\ge1$, of a generalized Bernstein
inequality using positive definite functions was recently proposed by
V.P.~Zastavnyi and A.D.~Manov~\cite{ZasMan18}. On homogeneous spaces $M=G/H$,
order Bernstein inequalities were studied by I.Z.~Pesenson~\cite{Pes91}. The
equality $\mathcal{C}_{pp}^{\partial^{r}}(\mathcal{T}_{n})=n^{r}$ for
fractional $r\ge1$ and $1\le p\le\infty$ was established by A.I.~Kozko (1998;
see the survey~\cite{AreGla14}). Numerous papers deal with variants of these
inequalities for polynomials and general multiplier differential operators
(see, for example,~\cite{Akh65,Are81,Bab92,MilMitRas94,RunSch01}).

Writing the $L^{\infty}$ Bernstein inequality for a polynomial
$P_{n}(\cos t)$, where $P_{n}\in\mathcal{P}_{n}$, gives the following
Markov-type inequality in $L^{\infty}([-1,1])$:
\[
|P_{n}'(x)|\le \frac{n}{\sqrt{1-x^{2}}}\,\|P_{n}\|_{\infty},
\quad x\in(-1,1),
\]
and for $x=0$ compare with~\eqref{Pnk-0}. The best currently known result for
the algebraic Markov inequality in $L^{p}([-1,1])$, $p\ge1$, has the form
\[
\mathcal{C}_{pp}^{\partial}(\mathcal{P}_{n})\le A(n,p)n^{2},
\]
with a nonsharp constant $A(n,p)$ (see~\cite{MilMitRas94}). Order
Bernstein--Nikolskii inequalities on different intervals were studied by
N.K.~Bari~\cite{Bar54}. The most general order results for the Markov brothers'
problem when $p\ge1$, including the case of a Gegenbauer algebraic weight, were
obtained by S.V.~Konyagin~\cite{Kon78}.

The next important step in the development of trigonometric inequalities in
approximation theory was made by D.~Jackson in 1933~\cite{Jac33}. Using
Bernstein's inequality, he derived the following Nikolskii inequality in
$L^{p}(\mathbb{T})$ for $p\ge1$:
\begin{equation}\label{jac-ineq}
\|T_{n}\|_{\infty}\le 2n^{1/p}\|T_{n}\|_{p},\quad n\ge1,
\end{equation}
and therefore
\begin{equation}\label{Cpinf-Tn}
\mathcal{C}_{p\infty}(\mathcal{T}_{n})\le 2n^{1/p}.
\end{equation}
Because of this result, $p\infty$ Nikolskii inequalities are also called
Jackson--Nikolskii inequalities.

Jackson's argument is as follows. Let $T\in\mathcal{T}_{n}$ and
$\|T\|_{\infty}=|T(t_{0})|$. For
$|t-t_{0}|\le\frac{1}{2n}$, the mean value theorem gives
\[
|T(t_{0})|-|T(t)|\le |T(t)-T(t_{0})|
\le \frac{\|T'\|_{\infty}}{2n}
\le \frac{\|T\|_{\infty}}{2}.
\]
Hence $|T(t)|\ge\frac{\|T\|_{\infty}}{2}$ on an interval of length at least
$\frac1n$. Therefore
$\|T\|_{p}^{p}\ge\frac{(\|T\|_{\infty}/2)^{p}}{n}$, which implies
\eqref{jac-ineq}.

Using the same idea together with the $L^{\infty}$ Markov inequality for
algebraic polynomials, Jackson proved that in $L^{p}([-1,1])$,
\[
\mathcal{C}_{p\infty}(\mathcal{P}_{n})
\le 2^{1+1/p}n^{2/p},\quad n\ge1,\quad p\ge1.
\]
He used these results to estimate uniform approximation of a continuous
function by best-approximation polynomials in other $L^{p}$ metrics (see
Subsection~\ref{subsec-appl-jac}).

For $p=1$,~\eqref{Cpinf-Tn} gives
\[
c_{n}=\mathcal{C}_{1\infty}(\mathcal{T}_{n})\le 2n,
\]
but, as we shall see below, this estimate is sharp only in order. We emphasize
this problem because it and related questions are now actively studied in
applications (see Section~\ref{sec-appl}).

In 1938, Ya.L.~Geronimus~\cite{Ger38} obtained the following remarkable result:
$c_{n}$ is the largest positive root of
\[
\begin{vmatrix}
1&0&\ldots&0&\mu_{0}&\mu_{1}&\ldots&\mu_{n}\\
0&1&\ldots&0&0&\mu_{0}&\ldots&\mu_{n-1}\\
\hdotsfor[2]{8}\\
0&0&\ldots&1&0&0&\ldots&\mu_{0}\\
\overline{\mu_{0}}&0&\ldots&0&1&0&\ldots&0\\
\overline{\mu_{1}}&\overline{\mu_{0}}&\ldots&0&0&1&\ldots&0\\
\hdotsfor[2]{8}\\
\overline{\mu_{n}}&\overline{\mu_{n-1}}&\ldots&\overline{\mu_{0}}&0&0&\ldots&1
\end{vmatrix}=0,
\]
where $\mu_{k}=\mu_{k}(c_{n})$ are defined by
\[
\tan
\Bigl\{\frac{1}{4c_{n}}\Bigl(\frac{1}{2}+\sum_{k=1}^{n}z^{k}\Bigr)\Bigr\}
=\mu_{0}+\mu_{1}z+\ldots+\mu_{n}z^{n}+\ldots\,.
\]
This is a beautiful result, but unfortunately it is not clear how to extract
from it simple estimates for $c_{n}$ in terms of $n$, including its asymptotic
behavior as $n\to\infty$. Geronimus's proof made essential use of complex
analysis. He also observed that the problem of finding $c_{n}$ is a
generalization of Chebyshev problems. The Geronimus results were extended to a
Markov inequality with Chebyshev weight by I.E.~Simonov and
P.Yu.~Glazyrina~\cite{SimGla15}.

In 1961 S.B.~Stechkin showed that there is a constant $c>0$ such that
\begin{equation}\label{st-cc}
c_{n}=cn+o(n),\quad n\to\infty.
\end{equation}
This result was not published by Stechkin himself and was stated without proof
in a paper by L.V.~Taikov~\cite{Tai65}, where the bounds
\begin{equation}\label{tai-bounds}
\frac{n}{\pi \operatorname{Si}\pi}+O(1)\le c_{n}\le
\frac{2Gn}{\pi^{2}}+O(1),\quad
\operatorname{Si}\pi=\int_{0}^{\pi}\frac{\sin t}{t}\,dt,
\end{equation}
were obtained, with $G$ denoting Catalan's constant. Hence
$c\in(0.17187,0.18562)$. The lower estimate was obtained using Rogozinski
polynomials. For the upper estimate, Taikov used the Riesz interpolation
formula in a duality relation going back to S.N.~Bernstein and closely related
to Bohr's inequality. Below we shall see how duality helps in more general
problems (see Subsection~\ref{subsec-common-dual}). Taikov later presented
Stechkin's proof of~\eqref{st-cc} in~\cite{Tai93}, a paper devoted to
approximating Dirichlet kernels by functions orthogonal to trigonometric
polynomials.

We now pause the discussion of $c_{n}$ and turn to entire functions of
exponential type. This subject began to develop actively around the 1930s
through the work of S.N.~Bernstein, M.G.~Krein, N.I.~Akhiezer, R.~Paley,
N.~Wiener, R.P.~Boas, B.Ya.~Levin, G.~P\'olya, G.~Szeg\H{o}, and many others.
Much of the theory developed in parallel with the corresponding problems for
trigonometric polynomials. Many classical results can be found in
\cite{Boa54,Tim60,Akh65,Lev96}.

Let $\mathcal{E}_{\sigma}=\mathcal{E}_{\sigma}^{\infty}$ be the set of all
entire functions of exponential type at most $\sigma>0$ that are bounded on
$\mathbb{R}$. Such functions $f(z)$ satisfy
$|f(z)|\le\|f\|_{\infty;\mathbb{R}}e^{\sigma|\!\operatorname{Im}z|}$ for every
$z\in\mathbb{C}$ (see~\cite{Lev96}). The least possible $\sigma$ is called the
type of $f$. The set $\mathcal{E}_{\sigma}$ is often called a Bernstein class.
In real analysis and approximation theory, an important role is played by the
classes of functions with finite $L^{p}$ norm on the real axis,
$\mathcal{E}_{\sigma}^{p}=\mathcal{E}_{\sigma}\cap L^{p}(\mathbb{R})$. For
every $p>0$ these classes are contained in $\mathcal{E}_{\sigma}$, since
$\|f\|_{\infty}\le C\|f\|_{p}$, where $C$ depends only on $p$ and $\sigma$
(a consequence of the Plancherel--P\'olya theorem~\cite{Lev96}). The Nikolskii
inequality is precisely a refinement of this constant.

The following analogue of the Riesz interpolation formula~\eqref{Tn-sum} holds
(see~\cite{Akh65}): for every $f\in\mathcal{E}_{\sigma}$,
\[
f'(x)=\sigma \sum_{k=-\infty}^{\infty}
\frac{(-1)^{k+1}}{(\pi k-\pi/2)^{2}}\,
f\Bigl(x+\frac{\pi k-\pi/2}{\sigma}\Bigr),\quad x\in\mathbb{R}.
\]
Applying H\"older's inequality, translation invariance of the
$L^{p}(\mathbb{R})$ norm, and
$\sum_{k=-\infty}^{\infty}(\pi k-\pi/2)^{-2}=1$, one obtains the sharp
Bernstein inequality in $L^{p}(\mathbb{R})$, $p\in[1,\infty]$, for derivatives
of entire functions of exponential type:
\begin{equation}\label{ber-ineq-ent}
\|f^{(r)}\|_{p}\le \sigma^{r}\|f\|_{p},\quad r\in\mathbb{N}.
\end{equation}
The difficult range $p\in(0,1)$ was treated by Q.I.~Rahman and
G.~Schmeisser~\cite{RahSch90}. For fractional $r\ge1$ and the Liouville
derivative, an inequality of the form~\eqref{ber-ineq-ent} for
$1\le p\le\infty$ was proved by P.I.~Lizorkin~\cite{Liz65}.

An analogue of Jackson's trigonometric inequality~\eqref{jac-ineq} for entire
functions of exponential type was established by J.~Korevaar in
1949~\cite{Kor49}:
\begin{equation}\label{kor-c}
\|f\|_{\infty}\le (A_{p}p\sigma)^{1/p}\|f\|_{p},\quad
f\in\mathcal{E}_{\sigma}^{p},
\end{equation}
where $A_{1}=\frac{1}{\pi}$,
$A_{p}=\frac{2^{k}}{p\pi}<\frac{1}{\pi}$,
$2^{k}<p\le2^{k+1}$, $k=0,1,2,\ldots.$ For the best constant $A_{p}^{*}$,
the bounds $\frac{1}{2p\pi}\le A_{p}^{*}\le\frac{1}{p\pi}$ were obtained for
$1\le p<2$, while $\frac{1}{p\pi}<A_{p}^{*}<\frac{1}{\pi}$ for $p>2$, and
$A_{2}^{*}=\frac{1}{2\pi}$ is the exact value.

\subsection{The era of Nikolskii}\label{subsec-hist-nik}
An important step in the theory of inequalities for polynomials and entire
functions of exponential type was made by S.M.~Nikolskii in 1951~\cite{Nik51}
(see also his book~\cite{Nik75}). It was largely after his work that the term
``Nikolskii inequalities'' came into use. Nikolskii considered multidimensional
versions of the Jackson--Korevaar inequality and obtained upper estimates not
only for the supremum norm but also for arbitrary $L^{q}$ norms.

Let $d\in\mathbb{N}$. For $x,y\in\mathbb{R}^{d}$, write
$xy=x_{1}y_{1}+\ldots+x_{d}y_{d}$ for the scalar product and
$|x|=(xx)^{1/2}$ for the Euclidean norm. Let
$\Omega\subset\mathbb{R}^{d}$ be a compact set (a body), $n\ge0$, and let
$\mathcal{T}_{n}(\Omega)$ be the set of trigonometric polynomials in $d$
variables
$T(x)=\sum_{k\in\mathbb{Z}^{d}\cap n\Omega}c_{k}e^{ikx}$ with complex
coefficients and spectrum contained in $n\Omega$. We restrict ourselves to
convex origin-symmetric bodies $\Omega$, with the cube
$\Omega=[-1,1]^{d}$ and the unit Euclidean ball
$B^{d}=\{x\in\mathbb{R}^{d}\colon |x|\le1\}$ as the main examples. The case of
nonconvex sets, such as hyperbolic crosses, is no less important and still
contains many open questions. Bernstein--Nikolskii inequalities and their
applications for hyperbolic crosses were studied by K.I.~Babenko,
S.A.~Telyakovskii, B.S.~Mityagin, and especially V.N.~Temlyakov (see
\cite{Tem86,TemTik17} and the references therein).

Classes of entire functions of exponential type with spectrum in
$\sigma\Omega$, $\sigma>0$, are defined in a similar way (see, for example,
\cite{SteWei74,Nik75,Gen77,NesWil78}; for a weighted version see
\cite{GorIvaTik19}). Namely,
$\mathcal{E}_{\sigma}(\Omega)=\mathcal{E}_{\sigma}^{\infty}(\Omega)$ denotes
the set of entire functions $f$ of $d$ complex variables whose restrictions to
$\mathbb{R}^{d}$ are bounded and whose $d$-dimensional Fourier transform
\[
\widehat{f}(y)=\frac{1}{(2\pi)^{d}}\int_{\mathbb{R}^{d}}f(x)e^{-ixy}\,dx,
\quad y\in\mathbb{R}^{d},
\]
understood here in the sense of distributions, is supported in
$\sigma\Omega$. Such functions satisfy
$f(z)=O(e^{\sigma\|\!\operatorname{Im}z\|_{\Omega^{*}}})$ for
$z\in\mathbb{C}^{d}$, where
$\|z\|_{\Omega^{*}}=\sup\{|zx|\colon x\in\Omega\}$ and $\Omega^{*}$ is the
polar body of $\Omega$. As in one dimension, for $p>0$ the class
$\mathcal{E}_{\sigma}(\Omega)$ contains the subclasses
$\mathcal{E}_{\sigma}^{p}(\Omega)=\mathcal{E}_{\sigma}(\Omega)\cap
L^{p}(\mathbb{R}^{d})$. Moreover, $\|f\|_{\infty}\le C\|f\|_{p}$, where $C$
depends on $d$, $p$, and $\sigma\Omega$. The sharp value of this constant is
given by the Nikolskii inequality. When $\Omega=B^{d}$, following Nikolskii we
speak of entire functions of spherical exponential type.

Nikolskii proved that for $1\le p<q\le\infty$,
\begin{equation}\label{nik-c-c}
\begin{aligned}
\mathcal{C}_{pq}(\mathcal{T}_{n}([-1,1]^{d}))&\le
(2n)^{d(1/p-1/q)},\\
\mathcal{C}_{pq}(\mathcal{E}_{\sigma}([-1,1]^{d}))&\le
(2\sigma)^{d(1/p-1/q)}.
\end{aligned}
\end{equation}
To prove these results he derived multidimensional interpolation formulas of
Riesz type. Their main applications were in approximation theory, including
inverse Jackson theorems and descriptions of Nikolskii--Besov classes (see
Subsection~\ref{subsec-appl-approx} and~\cite{Nik75}). Nikolskii considered
problems for parallelepipeds, but an affine change of variables reduces them to
the cube. The same applies to the ball and a general ellipsoid. We shall not
focus on this point here.

For compact $M$ and $q\le p$, H\"older's inequality applies, so the interesting
range for Nikolskii constants is $q>p$. For noncompact $M$, the condition
$q>p$ is also essential. Indeed, when $q=p$ the constant is $1$, while for
$q<p$ one may use the normalized Bessel function
$j_{\alpha}(x)=\frac{2^{\alpha}\Gamma(\alpha+1)J_{\alpha}(x)}{x^{\alpha}}$ on
$M=\mathbb{R}$, $\alpha\ge-1/2$ (its relevant properties are summarized, for
example, in~\cite{GorIva15}). The function $j_{\alpha}$ is an even entire
function of exponential type $1$, and
$j_{\alpha}(x)=C_{1}x^{-\alpha-1/2}(\cos(x-C_{2})+O(x^{-1}))$ as
$x\to\infty$. Thus, when $q<p$, one can choose $\alpha$ so that
$j_{\alpha}\in L^{p}(\mathbb{R})$ but
$j_{\alpha}\notin L^{q}(\mathbb{R})$. For Bernstein inequalities, however,
the range $q\le p$ is of interest. The case $q=p$ is classical; for $q<p$ and
trigonometric polynomials, see for example~\cite{Tai65-trud}.

Nikolskii also showed that the second inequality in~\eqref{nik-c-c} is sharp in
order with respect to~$\sigma$. For the lower estimate he used the Fej\'er
kernel
$f_{\sigma}(x)=\prod_{i=1}^{d}\bigl(\frac{\sin(\sigma x_{i}/2)}{x_{i}}\bigr)^2$,
for which a change of variables gives
$\frac{\|f_{\sigma}\|_{q}}{\|f_{\sigma}\|_{p}}=
\sigma^{d(1/p-1/q)}\frac{\|f_{1}\|_{q}}{\|f_{1}\|_{p}}$. In fact, if $D$ is
a homogeneous differential operator of order $r\ge0$, then for entire
functions of exponential type the functional
$\frac{\|Df\|_{q}}{\|f\|_{p}}$ is homogeneous of degree
$\lambda^{r+d(1/p-1/q)}$ under the dilation $f(\lambda x)$. Thus one may use
any nonzero function instead of the Fej\'er kernel, and
\[
\mathcal{C}_{pq}^{D}(\mathcal{E}_{\sigma}(\Omega))=
\sigma^{r+d(1/p-1/q)}
\mathcal{C}_{pq}^{D}(\mathcal{E}_{1}(\Omega)).
\]
In other words, in this setting it is enough to study the
Bernstein--Nikolskii inequality for $\sigma=1$. To show that the order of the
trigonometric constant is sharp, one usually uses the periodic Fej\'er kernel
$\prod_{i=1}^{d}\bigl(\frac{\sin(nx_{i}/2)}{\sin(x_{i}/2)}\bigr)^2$ (see
Section~\ref{sec-common}). The use of interpolation formulas to prove
multidimensional Bernstein--Nikolskii inequalities was further developed by
I.~Pesenson~\cite{Pes08,Pes09}.

The Nikolskii inequalities~\eqref{nik-c-c} were sharpened by I.I.~Ibragimov
and A.S.~Dzhafarov~\cite{IbrDzh61}:
\begin{equation}\label{ij-nik}
\mathcal{C}_{pq}(\mathcal{E}_{\sigma}([-1,1]^{d}))\le
\Bigl(\frac{p_{0}\sigma}{\pi}\Bigr)^{d(1/p-1/q)},\quad
p_{0}=\lceil \tfrac{p}{2}\rceil.
\end{equation}
For $d=1$ this estimate, together with several other results on
Bernstein--Nikolskii inequalities, appears in A.F.~Timan's well-known book on
approximation theory~\cite{Tim60}. It overlaps with Korevaar's
estimate~\eqref{kor-c}. The paper~\cite{IbrDzh61} also contains some
multidimensional Bernstein--Nikolskii inequalities for entire functions of
exponential type when $1\le p\le2$. We present the method behind them in a
general form in Subsection~\ref{subsec-common-repr}. For now, note that order
Bernstein--Nikolskii inequalities for all $p\ne q$ can be deduced from
$\mathcal{C}_{pq}^{D}(Y)\le
\mathcal{C}_{qq}^{D}(Y)\mathcal{C}_{pq}(Y)$ (see
Subsection~\ref{subsec-common-reduct}). For example, let
$Y=\mathcal{E}_{\sigma}([-1,1]^{d})$ and
$D=\partial_{x_{1}}^{r}$. The one-dimensional Bernstein inequality gives
$\mathcal{C}_{qq}^{D}(Y)\le\sigma^{r}$. Combining this with~\eqref{ij-nik}, we
obtain the order-sharp estimate
$\mathcal{C}_{pq}^{D}(Y)\le C\sigma^{r+d(1/p-1/q)}$ for $\sigma>0$, in
agreement with~\eqref{C-asymp}. Unfortunately, this argument does not yield
sharp constants, because the sharp Bernstein and Nikolskii inequalities have
different extremal functions.

For trigonometric polynomials, the Nikolskii inequality was sharpened in one
dimension by A.F.~Timan~\cite{Tim60} and Z.~Ziegler~\cite{Zie77}, and for
$d\ge1$ by I.I.~Ibragimov~\cite{Ibr58} and N.M.~Sabziev (1965; see
\cite{NesWil78}):
\begin{equation}\label{C-T-d}
\mathcal{C}_{pq}(\mathcal{T}_{n}([-1,1]^{d}))\le
\Bigl(\frac{2p_{0}n+1}{2\pi}\Bigr)^{d(1/p-1/q)}.
\end{equation}
This inequality is sharp for $p=2$, when $p_{0}=1$.

Historically,~\cite{Zie77} appeared considerably later than the works of
I.I.~Ibragimov and A.F.~Timan and treats only the case $q=\infty$. Z.~Ziegler,
however, gave a characterization of the extremal trigonometric polynomial for
$d=1$. Namely, for $1\le p<\infty$, an extremal polynomial $T_{n}^{*}$ for
$\mathcal{C}_{p\infty}(\mathcal{T}_{n})$ exists and is unique up to a constant
factor. It is even, all its $2n$ zeros lie in $(-\pi,\pi)$, and it satisfies
the orthogonality relations
\begin{equation}\label{Tn-orth}
\int_{\mathbb{T}}|T_{n}^{*}(t)|^{p-1}\operatorname{sign}T_{n}^{*}(t)
(\cos kt-1)\,dt=0,\quad k=1,2,\ldots,n.
\end{equation}
These relations also follow from general facts of approximation theory that
were published at about the same time in~\cite{Sha71,Kor76} (see also
Subsection~\ref{subsec-common-orth}). Ziegler also proved an interesting
interlacing property for the zeros of $T_{n}^{*}$ and $T_{n+1}^{*}$, which is
reminiscent of orthogonal polynomials~\cite{Seg62}. This nontrivial result was
substantially generalized in~\cite{PinZie79}. The constant
$\mathcal{C}_{2\infty}(\mathcal{T}_{n})$ was also computed in~\cite{Zie77}, but
this is already covered by~\eqref{C-T-d}.

For $p=1$, the characterization of the extremal polynomial had already been
used by Ya.L.~Geronimus~\cite{Ger38}. Relations of the form~\eqref{Tn-orth}
can be used to compute Nikolskii constants numerically for a given $n$, by
reducing~\eqref{Tn-orth} to a nonlinear system for the zeros of the extremal
polynomial. This may seem like a special question, but in view of the
connection between constants for polynomials and entire functions, and of the
applications of Nikolskii constants, it becomes more significant (see
Subsection~\ref{subsec-hist-now} and Section~\ref{sec-appl}). This approach was
developed in a series of papers by the author and I.A.~Martyanov
\cite{GorMar18,GorMar19,Mar20,GorMar20-cheb,GorMar20-cheb-alg,GorMar20-imm}.

One year earlier, in 1976, D.~Amir and Z.~Ziegler~\cite{AmiZie76} proved
related results for the one-dimensional Nikolskii constant for algebraic
polynomials, including orthogonality relations for $p\ge1$ and estimates that
we shall later call basic (see Subsection~\ref{subsec-common-repr}). In
particular, for $p=1$ and $n\in\mathbb{Z}_{+}$ in $L^{1}([-1,1])$,
\begin{equation}\label{az-bounds}
0.5(n+1)^{2}\ge M_{n}\ge 0.125\begin{cases}
(n+2)^{2},&\text{$n$ even},\\ (n+1)(n+3),&\text{$n$ odd},
\end{cases}
\end{equation}
where $M_{n}=\mathcal{C}_{1\infty}(\mathcal{P}_{n})$. A similar upper estimate
was obtained by T.K.~Ho (also in 1976; see~\cite{MilMitRas94}). These estimates
were improved in~\cite{GorMar20-cheb-alg} (see~\eqref{gm-nik-m}).

A fairly complete form of the estimates for Nikolskii constants for different
sets $\Omega$, covering~\eqref{nik-c-c} and~\eqref{ij-nik}, was proved in 1978
by R.J.~Nessel and G.~Wilmes in the well-known paper~\cite{NesWil78}, which
also contains a useful survey of earlier results. Unlike the preceding work,
it already includes the range where $p$ and $q$ are below $1$. If
$0<p\le q\le\infty$, then
\begin{equation}\label{nw-nik}
\begin{aligned}
\mathcal{C}_{pq}(\mathcal{T}_{n}(\Omega))&\le
\Bigl(\frac{|p_{0}n\Omega\cap\mathbb{Z}^{d}|}{(2\pi)^{d}}\Bigr)^{1/p-1/q},\\
\mathcal{C}_{pq}(\mathcal{E}_{\sigma}(\Omega))&\le
\Bigl(\frac{|p_{0}\sigma\Omega|}{(2\pi)^{d}}\Bigr)^{1/p-1/q}.
\end{aligned}
\end{equation}
Here $|A|$ denotes either the cardinality of a finite set $A$ or the measure of
a measurable set $A$; the relevant meaning is clear from the context. For
example, $|p_{0}\sigma[-1,1]^{d}|=(2p_{0}\sigma)^{d}$ (volume), which gives
\eqref{ij-nik}. The inequalities~\eqref{nw-nik} are sharp only for
$(p,q)=(2,\infty)$. A new and interesting case here is the ball
$\Omega=B^{d}$. Nikolskii, for example, studied this case by enclosing the ball
in cubes, which naturally causes a loss in the constants as the dimension $d$
grows.

The paper~\cite{NesWil78} also outlines a universal method for obtaining basic
estimates such as~\eqref{ij-nik} and~\eqref{nw-nik} from estimates of the
reproducing kernel of the subspace $Y$ (the Dirichlet kernel). We present this
method in Subsection~\ref{subsec-common-repr}. Substantial improvements of such
basic estimates with respect to different parameters, for example an
exponential improvement in the dimension $d$, are highly nontrivial. One
example is given below; see~\eqref{LA-d}. In this connection, recall that for
entire functions of exponential type, homogeneity of the functional defining
the Nikolskii constant allows us to restrict attention to $\sigma=1$.

One-dimensional trigonometric polynomials can equivalently be viewed on the
torus $\mathbb{T}=\mathbb{R}/2\pi\mathbb{Z}$ or on the unit circle
$\mathbb{S}^{1}$. In higher dimensions these lead to different constructions:
one may consider polynomials on the multidimensional torus
$\mathbb{T}^{d}=\mathbb{R}^{d}/2\pi\mathbb{Z}^{d}=(-\pi,\pi]^{d}$
(trigonometric polynomials) or on the multidimensional unit sphere
$\mathbb{S}^{d}=\{x\in\mathbb{R}^{d+1}\colon |x|=1\}$ (spherical polynomials).
Thus the structure of the space $M$ and harmonic analysis on it begin to play a
role. In the theory of Riemannian symmetric spaces, the torus
$\mathbb{T}^{d}$ is a compact space of rank~$d$, whereas the sphere
$\mathbb{S}^{d}$ belongs to the class of compact rank-one spaces, which also
includes projective spaces (see, for example,~\cite{Hel64,Lev98,ConSlo99}).

A similar dual viewpoint in harmonic analysis appears for the Euclidean space
$\mathbb{R}^{d}$, $d\ge2$: it may be viewed either as an abelian group under
addition or as the quotient of the Euclidean motion group by the subgroup of
proper rotations $SO(d)$.

Order estimates for Nikolskii constants on the sphere $\mathbb{S}^{d}$ were
obtained by A.I.~Kamzolov~\cite{Kam84}. Let $\Pi_{n}^{d}$ be the subspace of
spherical polynomials of order at most~$n$, that is, restrictions to the
sphere of algebraic polynomials
\begin{equation}\label{poly-d}
f(x_{1},\ldots,x_{d+1})=\sum_{\substack{k_{1},\ldots,k_{d+1}\in
\mathbb{Z}_{+}\\ k_{1}+\ldots+k_{d+1}\le n}}c_{k_{1}\ldots
k_{d+1}}x_{1}^{k_{1}}\ldots x_{d+1}^{k_{d+1}}
\end{equation}
of degree at most $n$. Then, for fixed $d$ and $p,q\ge1$,
\[
\mathcal{C}_{pq}(\Pi_{n}^{d})\asymp n^{d(1/p-1/q)_{+}},\quad n\to\infty.
\]
Many order versions of classical inequalities on the multidimensional sphere,
including weighted versions, were later obtained by F.~Dai and Yu.~Xu
\cite{DaiXu13}. Kamzolov also treated spherical harmonics of order $n$, a
particularly difficult case. The harmonic subspace was studied in
\cite{DaiFenTik16}; see also the references there. We also mention
\cite{DaiGorTik20}, where Nikolskii inequalities are proved for lacunary
polynomials on the sphere, a class that in a certain sense lies between the
harmonics of order~$n$ and all spherical polynomials of order at most~$n$.
Another relevant paper is that of M.V.~Deikalova~\cite{Dei09}, who obtained an
upper estimate for $\mathcal{C}_{1\infty}(\Pi_{n}^{d})$ at $p=1$ equivalent to
the reproducing-kernel estimate from Subsection~\ref{subsec-common-repr}. The
case of general $p$ can be deduced from~\cite{BelDaiDit03}. We present this
approach in Subsection~\ref{subsec-common-repr}.

Earlier, in 1974, A.I.~Kamzolov~\cite{Kam74} studied Bernstein inequalities on
$d$-dimensional compact rank-one manifolds $M=M_{d}$. Let
$\Pi_{n}(M)=\bigoplus_{l=0}^{n}V_{l}$ be the subspace of polynomials of degree
at most $n$ on $M$, where $V_{l}$ are invariant subspaces (harmonics) carrying
irreducible class-one representations of the motion group of $M$. Let
$\Delta_{M}$ be the Laplace--Beltrami operator on $M$; its fractional powers on
$\Pi_{n}(M)$ are defined as multipliers by
$\Delta_{M}V_{l}=-\lambda_{l}V_{l}$, where
$\lambda_{0}=0<\lambda_{1}<\lambda_{2}<\ldots$ are the eigenvalues. We have
$\dim\Pi_{n}(M)\asymp n^{d}$ and $\lambda_{l}\asymp l^{2}$. Kamzolov proved an
analogue of the Riesz interpolation formula and deduced
$\mathcal{C}_{pp}^{\Delta_{M}}(\Pi_{n}(M))\le c_{M}dn^{2}$,
$1\le p\le\infty$, where $c_{M}=1,2,4$ for the sphere, the complex projective
space, and the quaternionic projective space, respectively. This estimate is
sharp for $p=\infty$. Kamzolov also considered the Euclidean case, which we
mention below.

In 1983 V.A.~Ivanov~\cite{Iva83} proved a general order
Bernstein--Nikolskii inequality for compact rank-one manifolds when
$p,q\ge1$: for $r\ge0$ and $n\to\infty$,
\[
\mathcal{C}_{pq}^{(-\Delta_{M})^{r/2}}(\Pi_{n}(M))
\asymp n^{r+d(1/p-1/q)_{+}}.
\]
For the sphere and $p=q$, this result was also obtained by
A.I.~Kamzolov~\cite{Kam84-umn}. The result was later extended to the sphere,
including weighted settings, and to $p>0$ by F.~Dai and S.~Tikhonov
\cite{Dai06,DaiTik16}.

In 1992 V.A.~Ivanov~\cite{Iva92} computed the constants
$\mathcal{C}_{pq}(\Pi_{n}(M),(-\Delta_{M})^{r/2})$ for
$(p,q)=(2,2)$, $(\infty,\infty)$, and $(2,\infty)$. For the cases
$(2,2)$ and $(2,\infty)$ there is a general method, described in
Subsection~\ref{subsec-common-repr}, so we do not reproduce the values here. He
also computed the Bernstein constant
$\mathcal{C}_{\infty\infty}^{(-\Delta_{M})^{r/2}}(\Pi_{n}(M))$ for every
integer $r$ and for fractional $r\ge r_{0}$; it was conjectured that
$r_{0}=1$. The case $r=2$ had been obtained in~\cite{Kam74}. Unlike the Riesz
formula approach of~\cite{Kam74}, Ivanov used, for integer $r$, an original
argument reducing the problem to the general Markov brothers'
inequality~\eqref{mar-ineq} and the remarkable fact that the same Chebyshev
polynomial of the first kind is extremal for all derivatives. The explicit
answer is rather cumbersome; the main point is that the extremal polynomial is
the zonal polynomial $\cos n\rho(x,o)$, where $\rho$ is the metric on $M$.

As Ivanov notes~\cite{Iva92}, in one dimension and for fractional $r$ one
obtains the Riesz derivative
$(-\partial_{t})^{r/2}T_{n}(t)=\sum_{k=-n}^{n}|k|^{r}c_{k}e^{ikt}$, which
differs from the fractional Weyl derivative
$T_{n}^{(r)}(t)=\sum_{k=-n}^{n}(ik)^{r}c_{k}e^{ikt}$; they agree for even $r$.
For the Weyl derivative, the sharp Bernstein inequality
$\|T_{n}^{(r)}\|_{p}\le n^{r}\|T_{n}\|_{p}$ for all
$1\le p\le\infty$ and $r\ge1$ was proved by A.I.~Kozko, as mentioned above.
For $0<r<1$ only estimates of the best constant are known. The case $p<1$ for
noninteger $r$ remains open; see the survey~\cite{AreGla14}.

The ideas of~\cite{Iva92} were used in~\cite{GorIva19} to compute the sharp
Bernstein constant in $L^{\infty}(\mathbb{R}^{d})$ for entire functions of
spherical exponential type and powers of the differential-difference
Dunkl--Laplace operator. For the first power of the ordinary Laplacian and for
spaces of entire functions of exponential type corresponding to the cube and
the ball, this had been done by A.I.~Kamzolov in 1974~\cite{Kam74}. These
results follow from the one-dimensional Bernstein inequality, but Kamzolov also
gave an alternative proof based on an analogue of the Riesz interpolation
formula. We also mention the later paper of S.S.~Platonov~\cite{Pla07}, where a
Riesz-type interpolation formula is constructed on the half-line
$\mathbb{R}_{+}$ with the power weight $x^{2\alpha+1}$,
$\alpha\ge-1/2$. This gives sharp $pp$ Bernstein inequalities for
$p=2,\infty$. For other $p\ne2,\infty$, sharp multidimensional and weighted
Bernstein inequalities are still unknown. One-dimensional weighted order
Bernstein inequalities can be found in D.S.~Lubinsky's paper~\cite{Lub14}.
General order Bernstein--Nikolskii inequalities in
$L_{v}^{p}(\mathbb{R}^{d})$ with a power Dunkl weight $v$, for entire functions
of spherical exponential type and powers of the Dunkl--Laplace operator, were
proved recently by the author together with V.I.~Ivanov and S.Yu.~Tikhonov
\cite{GorIvaTik19,GorIva19-mz}.

\subsection{The modern era}\label{subsec-hist-now}
We now return to the trigonometric constant

$c_{n}=\mathcal{C}_{1\infty}(\mathcal{T}_{n})$. Let
$L=\mathcal{C}_{1\infty}(\mathcal{E}_{1})$ be the corresponding constant for
entire functions of exponential type at most $1$. In 2005 the author proved in
\cite{Gor05} that for every $n\in\mathbb{Z}_{+}$,
\begin{equation}\label{cL}
nL\le c_{n}\le(n+1)L.
\end{equation}
This immediately gives
\begin{equation}\label{dg-cL}
c_{n}=Ln(1+O(n^{-1})),\quad n\to\infty,
\end{equation}
which slightly strengthens Stechkin's result~\eqref{st-cc}. It was also shown
in~\cite{Gor05} that
$L=\sup\frac{\widehat{\varphi}(0)}{\|\widehat{\varphi}\|_{1}}$, where the
supremum is taken over all even functions $\varphi\in C(\mathbb{R})$ with
$\operatorname{supp}\varphi\subset[-1,1]$. In this form the problem was studied
in~\cite{AndKonPop96} (with the correction~\cite{AndKonPop00}), where the
bounds
$\frac{1}{\pi\operatorname{Si}\pi}\le L\le
\frac{\operatorname{Si}\pi}{\pi^{2}}$ were proved. Thus
$L\in(0.17187,0.18765)$; compare with Taikov's bounds~\eqref{tai-bounds}. In
\cite{Gor05} the lower and upper bounds were improved to
$L\in(0.17218,0.17471)$. The recent paper~\cite{CarMilSou19} emphasizes the
importance of computing the Nikolskii constant $L$ for applications in
analytic number theory; see Subsection~\ref{subsec-appl-num}.

A numerical value of $L$ to several decimal places was obtained in 1993 by
L.~H\"ormander and B.~Bernhardsson~\cite{HorBer93} in their derivation of a
two-dimensional Bohr inequality: $L\approx0.172182$. They intended to determine
$L$ exactly, but were unable to do so. The problem remains open.

A different estimate of $L$, based on~\eqref{cL} and~\eqref{Tn-orth}, was
given in~\cite{GorMar18}. Indeed,~\eqref{cL} gives
$\frac{c_{n}}{n+1}\le L\le\frac{c_{n}}{n}$. For sufficiently large $n$, the
constant $c_{n}$ can be computed approximately from~\eqref{Tn-orth} with
$p=1$, written as a system of nonlinear equations for the unknown zeros of the
extremal polynomial. Once these zeros are found, $c_{n}$ is recovered from
formulas similar to those in~\cite{AmiZie76,AreDei13}. This method was used in
\cite{GorMar20-cheb-alg,GorMar20-imm} to sharpen trigonometric and algebraic
Bernstein--Nikolskii constants. For example,
$M_{n}=\mathcal{C}_{1\infty}(\mathcal{P}_{n})=
2\mathcal{C}_{1\infty}^{\partial}(\mathcal{T}_{n+1})$, and for
$n\in\mathbb{Z}_{+}$,
\[
n^{2}\mathcal{C}_{1\infty}^{\partial}(\mathcal{E}_{1})\le
\mathcal{C}_{1\infty}^{\partial}(\mathcal{T}_{n})\le
(n+1)^{2}\mathcal{C}_{1\infty}^{\partial}(\mathcal{E}_{1}).
\]
This led in~\cite{GorMar20-cheb-alg} to the estimates
\begin{equation}\label{gm-nik-m}
(n+1)^{2}M\le M_{n}\le(n+2)^{2}M,\quad M\in(0.1440,0.1441),
\end{equation}
which are stronger than~\eqref{az-bounds}, where
$M=2\mathcal{C}_{1\infty}^{\partial}(\mathcal{E}_{1})$. As with $L$, the exact
value of $M$ is unknown. The asymptotic formula
$M_{n}=Mn^{2}(1+o(1))$ and the bound $M\in(0.141,0.192)$ were also proved in
\cite{DaiGorTik21}, which deals with the multidimensional sphere.

In our view, a major advance in sharp Bernstein--Nikolskii inequalities is the
series of results establishing asymptotic equalities of the form~\eqref{ll-cL}.
In 2015 E.~Levin and D.~Lubinsky~\cite{LevLub15} (see also
\cite{LevLub15-anal}) extended Stechkin's result~\eqref{st-cc} and the author's
formula~\eqref{dg-cL} to $L^{p}$ spaces for all $p>0$. They proved that
\begin{equation}\label{ll-eq}
c_{np}=L_{p}n^{1/p}(1+o(1)),\quad n\to\infty,\quad0<p<\infty,
\end{equation}
where $c_{np}=\mathcal{C}_{p\infty}(\mathcal{T}_{n})$ and
$L_{p}=\mathcal{C}_{p\infty}(\mathcal{E}_{1})$. In a sense, this result unifies
the problems for trigonometric polynomials and entire functions of exponential
type. For comparison, the only known exact formula corresponds to $p=2$:
$c_{n2}=(n+1/2)^{1/2}L_{2}$, where $L_{2}=1/\sqrt{\pi}$; compare with
Korevaar's estimate~\eqref{kor-c}. Technically, the proof of~\eqref{ll-eq}
uses periodization based on the Poisson summation formula.

We note at once that no result of the form~\eqref{ll-eq} is known for
$q\ne\infty$; only the required lower estimate is available. Levin and
Lubinsky explicitly pointed out this important open problem.

We also mention D.~Lubinsky's paper~\cite{Lub14-ams}, where a result similar in
spirit to~\eqref{ll-eq} is proved for the relation between constants in
Marcinkiewicz--Zygmund and Plancherel--P\'olya inequalities.

In analogy with~\eqref{cL}, the result was strengthened in~\cite{GorMar18}.
For every $p>0$ and $n\in\mathbb{Z}_{+}$,
\begin{equation}\label{cL-p}
n^{1/p}L_{p}\le c_{np}\le
(n+\lceil\tfrac1p\rceil)^{1/p}L_{p}.
\end{equation}
We shall prove these inequalities in even greater generality in
Section~\ref{sec-ll}. In this connection we mention the conceptually similar
inequality proved in~\cite{AshGan99}:
\[
\alpha_{p}\le A_{np}\le\alpha_{p}(1+n^{-1}),\quad1\le p\le\infty,
\]
where $\alpha_{p}=\pi\|\!\cos t\|_{p'}^{-1}$, $p'=p/(p-1)$ is the conjugate
exponent, and $A_{np}=\inf\|T_{n}\|_{p}$, with the infimum taken over all real
polynomials
$T_{n}(x)=\sum_{k=0}^{n}(a_{k}\cos kx+b_{k}\sin kx)$ satisfying
\[
\|\widehat{T}_{n}\|_{\infty}=
\max\{|a_{0}|,|a_{1}|,\ldots,|a_{n}|,|b_{1}|,\ldots,|b_{n}|\}=1.
\]
Moreover, $A_{np}=\alpha_{p}$ if $p'=2,4,\ldots$ and $n\ge p'-1$. The
similarity lies in the fact that
$c_{np}^{-1}=\inf\|T_{n}\|_{p}$ under the normalization
$\|T_{n}\|_{\infty}=1$; the normalizations themselves are different.

The Levin--Lubinsky result was extended to Bernstein--Nikolskii inequalities by
M.~Ganzburg and S.~Tikhonov~\cite{GanTik17}: for every
$r\in\mathbb{Z}_{+}$,
\begin{equation}\label{gt-eq}
c_{np}^{(r)}=L_{p}^{(r)}n^{r+1/p}(1+o(1)),
\end{equation}
where
$c_{np}^{(r)}=\mathcal{C}_{p\infty}^{\partial^{r}}(\mathcal{T}_{n})$ and
$L_{p}^{(r)}=\mathcal{C}_{p\infty}^{\partial^{r}}(\mathcal{E}_{1})$. They also
proved the existence of extremal polynomials
$\tilde{T}_{nr}\in\mathcal{T}_{n}$ and extremal functions
$\tilde{F}_{r}\in\mathcal{E}_{1}^{p}$ with unit $p$-norms such that
$c_{np}^{(r)}=\tilde{T}_{nr}^{(r)}(0)$ and
$L_{p}^{(r)}=\tilde{F}_{r}^{(r)}(0)$. Existence of an extremizer can itself be
nontrivial; the general issue is discussed in
Subsection~\ref{subsec-common-exist}. The proof in~\cite{GanTik17} uses
properties of Levitan polynomials constructed from the Poisson summation
formula.

The Ganzburg--Tikhonov result~\eqref{gt-eq} was sharpened in~\cite{GorMar19}
and in the subsequent correction~\cite{GorMar20-cheb-lett}. Let, for
$s\in\mathbb{N}$,
\[
\Bigl(\frac{\sin\pi x}{\pi x}\Bigr)^{2s}
=\sum_{k=0}^{\infty}\frac{(-1)^{k}A_{sk}(2\pi x)^{2k}}{(2k)!},
\quad A_{s0}=1,\quad A_{sk}>0.
\]
Then, for $s=\lceil\frac{r+1}{2}\rceil$,
\[
c_{np}^{(r)}\ge(n-s+1)^{r+1/p}\biggl(L_{p}^{(r)}+
\sum_{i=1}^{\lfloor r/2\rfloor}
\frac{(-1)^{i}\binom{r}{2i}A_{si}}{n^{2i}}\,
\tilde{F}_{r}^{(r-2i)}(0)\biggr),
\]
and, for $s=\lceil\frac1p\rceil$,
\[
c_{np}^{(r)}+
\sum_{i=1}^{\lfloor r/2\rfloor}(-1)^{i}\binom{r}{2i}A_{si}
\tilde{T}_{nr}^{(r-2i)}(0)
\le(n+s)^{r+1/p}L_{p}^{(r)},
\]
where
$|\tilde{T}_{nr}^{(r-2i)}(0)|\le
n^{r-2i}(n+\lceil\frac1p\rceil)^{1/p}L_{p}$. These inequalities provide
estimates for the remainder terms in~\eqref{gt-eq}.

The papers~\cite{GorMar19,GorMar20-cheb-lett} propose the open conjecture that
the signs of the nonzero Taylor coefficients
$\tilde{T}_{nr}^{(j)}(0)$ and $\tilde{F}_{r}^{(j)}(0)$ alternate. If this is
true, then for every $p\in(0,\infty]$, $r\in\mathbb{Z}_{+}$, and
$n\ge\lceil\frac{r+1}{2}\rceil-1$ ($n\ge1$ when $r=1$),
\[
(n-\lceil\tfrac{r+1}{2}\rceil+1)^{r+1/p}L_{p}^{(r)}
\le c_{np}^{(r)}\le
(n+\lceil\tfrac1p\rceil)^{r+1/p}L_{p}^{(r)}.
\]

Analogous results for the algebraic Nikolskii constant follow from general
inequalities for the Gegenbauer weight~\cite{GorMar20-imm}; see~\eqref{A-nn}
and Section~\ref{sec-ll}.

The next important step toward Levin--Lubinsky type asymptotic formulas was made
in the author's joint paper~\cite{DaiGorTik20-anal}, where the multidimensional
sphere $\mathbb{S}^{d}$ and spherical polynomials were considered. For fixed
$d\in\mathbb{N}$ and every $0<p<\infty$,
\begin{equation}\label{ll-S}
\mathcal{C}_{p\infty}(\Pi_{n}^{d})=
\mathcal{L}_{dp}n^{d/p}(1+o(1)),\quad n\to\infty,
\end{equation}
where
$\mathcal{L}_{dp}=\mathcal{C}_{p\infty}(\mathcal{E}_{1}(B^{d}))$ is the
Nikolskii constant in $L^{p}(\mathbb{R}^{d})$ for entire functions of spherical
exponential type at most $1$.

There is no known spherical analogue of periodization by the Poisson summation
formula, so the proof of~\eqref{ll-S} is technically more difficult than on the
torus. It required nontrivial tools from spherical harmonic analysis, including
a special approximation to the reproducing kernel of the space of spherical
polynomials, the deep result of~\cite{BonRadVia15} on the existence of
well-distributed spherical designs, estimates for such distributions, and
Marcinkiewicz--Zygmund type results for spherical polynomials.

For multidimensional trigonometric and algebraic polynomials,
formula~\eqref{ll-S} was generalized by M.~Ganzburg in a series of papers
\cite{Gan20-jfaa,Gan20,Gan21-jmaa,Gan21-jat,Gan21}. He obtained many
extensions of Levin--Lubinsky formulas for multidimensional
Bernstein--Nikolskii problems, including shifted problems; see
Subsection~\ref{subsec-common-shift}. In particular, for trigonometric
polynomials with spectrum in a ball one has the analogue of~\eqref{ll-S}
\cite{Gan20-jfaa}
\begin{equation}\label{ll-T}
\mathcal{C}_{p\infty}(\mathcal{T}_{n}(B^{d}))=
\mathcal{L}_{dp}n^{d/p}(1+o(1)),\quad p>0.
\end{equation}
It follows from the more general equality
\[
\mathcal{C}_{p\infty}^{D_{s}}(\mathcal{T}_{n}(\Omega))=
\mathcal{C}_{p\infty}^{D_{s}}(\mathcal{E}_{1}(\Omega))
 n^{s+d/p}(1+o(1)),
\]
where $D_{s}$ is a homogeneous differential operator of order
$s\in\mathbb{Z}_{+}$, that is, a homogeneous polynomial in
$\partial_{x_{1}},\ldots,\partial_{x_{d}}$ with constant coefficients.

Analogous results for the sphere and spherical polynomials, and for the cube
and algebraic polynomials, can be found in~\cite{Gan19}, but only for
$p\ge1$, when the maximum can be shifted: in taking the supremum,
$\|Df\|_{\infty}$ can be replaced by the value $|Df(o)|$ at a fixed point; see
Subsection~\ref{subsec-common-shift}. For $p<1$, shifting the maximum becomes
problematic. Thus, for example, in~\cite{Gan21-jmaa} a Levin--Lubinsky type
formula in $L^{p}([-1,1]^{d})$ for algebraic polynomials is proved for $p>0$
only for a Markov coefficient problem (see~\eqref{Pnk-0}), which is naturally
formulated in terms of derivatives at the origin. Proving the main
multidimensional Markov--Nikolskii inequalities for $p<1$ remains open. It has
been solved only for $d=1$; see~\cite{GorMar20-imm}.

In view of~\eqref{ll-S} and~\eqref{ll-T}, good estimates for the Nikolskii
constant $\mathcal{L}_{dp}$ become especially important. The basic upper bound
follows from~\eqref{nw-nik} with $\sigma=1$ and $\Omega=B^{d}$:
\[
\mathcal{L}_{dp}\le
\Bigl(\frac{|B^{d}|}{(2\pi)^{d}}\Bigr)^{1/p}
\lceil\tfrac{p}{2}\rceil^{d/p},\quad p>0.
\]
For $p=2$ this estimate is sharp. The results of~\cite{GorDob18} imply that,
for fixed $p\ge1$ and $d\to\infty$,
\[
\mathcal{L}_{dp}\ge
\Bigl(\frac{|B^{d}|}{(2\pi)^{d}}\Bigr)^{1/p}
(\tfrac{p}{2})^{d(1+o(1))/p}.
\]
For even $p\ge2$, the upper and lower estimates match. The case $p=1$ is
particularly intriguing. Note that $\mathcal{L}_{11}=L$, so the bounds for $L$
can be used in dimension one. In general, for $d\ge1$,
\[
\mathcal{L}_{d1}\le A_{d}=\frac{|B^{d}|}{(2\pi)^{d}}.
\]
The results of~\cite{Dei09} imply
$\mathcal{L}_{d1}\ge e^{-d(1+o(1))}A_{d}$ as $d\to\infty$. Thus the known
bounds for $\mathcal{L}_{d1}/A_{d}$ differ exponentially. They were
substantially improved in~\cite{DaiGorTik21}, where it was proved that for every
$d\ge1$,
\begin{equation}\label{LA-d}
2^{-d}\le \frac{\mathcal{L}_{d1}}{A_{d}}\le
{}_{1}F_{2}\Bigl(\frac{d}{2};\frac{d}{2}+1,\frac{d}{2}+1;
-\frac{\beta_{d}^{2}}{4}\Bigr)
=(\sqrt{2/e})^{d(1+\varepsilon_{d})},
\end{equation}
where $\sqrt{2/e}=0.857\ldots$ and
$\varepsilon_{d}=O(d^{-2/3})$ as $d\to\infty$. For $d=1$ this gives
$\mathcal{L}_{11}\le0.18764\ldots$, an estimate already mentioned above.

To prove the upper estimate in~\eqref{LA-d}, the paper~\cite{DaiGorTik21}
solved an extremal problem of $L^{\infty}$ approximation of the reproducing
kernel of $\mathcal{E}_{1}^{1}(B^{d})$ by Fourier--Bessel expansions. This
problem comes from duality; see Subsection~\ref{subsec-common-dual}. The proof
of the lower estimate in~\eqref{LA-d} is based on solving a Tur\'an extremal
problem by means of a quadrature formula at the zeros of the reproducing kernel
(a Bessel function). For the polynomial Nikolskii constant
$\mathcal{C}_{p\infty}(\Pi_{n}^{d})$, a related lower estimate was proposed in
\cite{DaiGorTik20-anal}; it is also proved by quadrature formulas and, up to a
factor, coincides with the famous Delsarte--Goethals--Seidel bound for tight
designs; see Subsection~\ref{subsec-appl-geom}. The paper~\cite{DaiGorTik21}
also formulates dual Nikolskii problems for spherical polynomials and entire
functions of spherical exponential type when $1\le p<\infty$, and proves the
existence and uniqueness of a radial extremal function in the problem
$\mathcal{L}_{d1}$. General ideas behind such results are discussed in
Section~\ref{sec-common}.

Among other ingredients, M.~Ganzburg's general results use ideas developed by
V.V.~Arestov and M.V.~Deikalova in their work on the spherical Nikolskii
constant and related problems~\cite{AreDei13,AreDei15,AreDei16}. They appear to
have been among the first to use generalized translation operators to obtain a
shifted Nikolskii problem; see Subsection~\ref{subsec-common-shift}. In
\cite{AreDei13}, the spherical Nikolskii constant
$\mathcal{C}_{p\infty}(\Pi_{n}^{d})$ was reduced to a one-dimensional weighted
problem with a Gegenbauer algebraic weight, which turned out to be connected
with a weighted Chebyshev problem. To state these results, let
$L_{w}^{p}(M)$ be the weighted space with
$\|f\|_{p,w}=(\int_{M}|f|^{p}w\,dx)^{1/p}$ for $p<\infty$, where
$w\colon M\to\mathbb{R}_{+}$ is a weight. As before, for $p=\infty$ the weight
does not enter the norm. We denote the corresponding weighted
Bernstein--Nikolskii constant by $\mathcal{C}_{pq,w}^{D}(Y)$.

For $1\le p<\infty$,~\cite{AreDei13} gives
\[
|\mathbb{S}^{d-1}|^{1/p}\mathcal{C}_{p\infty}(\Pi_{n}^{d})=
\mathcal{C}_{p\infty,(1-x^{2})^{d/2-1};[-1,1]}(\mathcal{P}_{n})
=\sup_{P\in\mathcal{P}_{n}}
\frac{|P(1)|}{\|P\|_{p,(1-x^{2})^{d/2-1}}}.
\]
The last supremum is attained by the extremal polynomial
$P_{n}^{*}(x)=x^{n}+\ldots$, which is the unique extremal polynomial in the
Chebyshev problem for the Jacobi weight:
\[
P_{n}^{*}=\operatorname*{argmin}_{P\in\mathcal{P}_{n-1}}
\|x^{n}+P(x)\|_{p,w},\quad
w(x)=(1-x)^{d/2}(1+x)^{d/2-1}.
\]
Viewed as a zonal polynomial, $P_{n}^{*}$ is also the unique extremizer, up to a
constant factor and a rotation of the argument, for the Nikolskii constant
$\mathcal{C}_{p\infty}(\Pi_{n}^{d})$ when $1\le p<\infty$.

We began this section with Chebyshev's problem of a polynomial of least
deviation from zero, and have now returned to it in a weighted form. Exact
values of the Nikolskii constant, as well as the solution of the Chebyshev
problem for a general Jacobi weight, are known only for $p=2$. In that case,
as noted above, the solution of the weighted Chebyshev problem is the
orthogonal polynomial for the corresponding weight.

The connection with the classical Chebyshev problem allowed the authors of
\cite{AreDei13} to use standard facts from approximation theory concerning
best-approximation polynomials. In particular, they obtained the following
orthogonality relation characterizing the extremal polynomial:
\[
\int_{-1}^{1}P|P_{n}^{*}|^{p-1}\operatorname{sign}P_{n}^{*}\,w\,dt=0,
\quad \forall\,P\in\mathcal{P}_{n-1}.
\]
It follows, in particular, that all zeros of $P_{n}^{*}$ are simple and lie in
$(-1,1)$. These facts generalize the results of D.~Amir and
Z.~Ziegler~\cite{AmiZie76} for $d=1$.

In~\cite{AreDei15,AreDei16}, V.V.~Arestov and M.V.~Deikalova extended the
relation between weighted Nikolskii and Chebyshev problems on $[-1,1]$ to an
arbitrary Jacobi weight
$w_{\alpha,\beta}(x)=(1-x)^{\alpha}(1+x)^{\beta}$ with
$\alpha\ge\beta\ge-1/2$, including the Gegenbauer case $\alpha=\beta$. They
proved that for $1\le p<\infty$,
\begin{equation}\label{C-P-cheb}
\mathcal{C}_{p\infty,w_{\alpha,\beta}}(\mathcal{P}_{n})=
\sup_{P\in\mathcal{P}_{n}}
\frac{|P(1)|}{\|P\|_{p,w_{\alpha,\beta}}},
\end{equation}
and equality is attained by the unique extremal polynomial from the weighted
Chebyshev problem
$P_{n}^{*}=\operatorname*{argmin}_{P\in\mathcal{P}_{n-1}}
\|x^{n}+P(x)\|_{p,w_{\alpha+1,\beta}}$. A substantial difficulty, already
mentioned above, was to prove that in computing the $p\infty$ Nikolskii
constant for algebraic polynomials, the norm $\|P\|_{\infty}$ may be replaced
by the value $|P(1)|$ at a fixed point; this is the shift of the maximum. This
fact is needed, in particular, to characterize the extremal polynomial in
\eqref{C-P-cheb}. For homogeneous spaces $M=G/H$, it is easy to establish:
the motion group $G$ gives a natural translation $f(gx)$ under which the
measure, the $L^{p}$ norm, and the class $Y$ are invariant; see
Subsection~\ref{subsec-common-shift}. On $[-1,1]$ with a general Jacobi weight,
there is no natural translation. Arestov and Deikalova successfully replaced
it with a symmetric positive generalized Jacobi translation operator. This
approach was further developed in
\cite{AreBabDeiHor18,GorDob18,GorIva19,Mar20,Gan19,Gan21}; we describe the
main idea in Subsection~\ref{subsec-common-shift}. We especially mention the
case of nonsymmetric Dunkl-type generalized translations, where additional
difficulties arise. In particular,~\cite{GorDob18,Mar20} proved, for
$1\le p<\infty$ and $\alpha\ge-1/2$, using positive Gegenbauer--Dunkl and
Bessel--Dunkl operators, that for $M=\mathbb{T}$ and $M=\mathbb{R}$,
respectively,
\begin{equation}\label{C-0}
\begin{aligned}
\mathcal{C}_{p\infty,|\!\sin x|^{2\alpha+1}}(\mathcal{T}_{n})&=
\sup_{T\in\mathcal{T}_{n}^{\text{real},\text{even}}}
\frac{|T(0)|}{\|T\|_{p,|\!\sin x|^{2\alpha+1}}},\\
\mathcal{C}_{p\infty,|x|^{2\alpha+1}}(\mathcal{E}_{1})&=
\sup_{f\in\mathcal{E}_{1}^{\text{real},\text{even}}\cap
L_{|x|^{2\alpha+1}}^{p}(\mathbb{R})}
\frac{|f(0)|}{\|f\|_{p,|x|^{2\alpha+1}}}.
\end{aligned}
\end{equation}
After the trigonometric substitution, the first equality reduces to the
algebraic Nikolskii constant for the Gegenbauer weight
$(1-x^{2})^{\alpha}$. For Bernstein--Markov--Nikolskii constants on
$M=\mathbb{T}$ and $[-1,1]$, see~\cite{Gan19,GorMar20-cheb}. In
Subsection~\ref{subsec-common-shift}, we call the right-hand sides of
\eqref{C-0} shifted Nikolskii problems and use, for example, the notation
$\mathcal{A}_{p,|\!\sin x|^{2\alpha+1}}
(\mathcal{T}_{n}^{\text{real},\text{even}};0)$. Analogues for
$q\ne\infty$ are not known.

For $p\ge1$, shifting the maximum in $p\infty$
Bernstein--Markov--Nikolskii constants resolves several auxiliary issues: it
reduces the problem to a primal-dual convex real problem, gives orthogonality
relations, characterizes the extremal function, and consequently proves its
uniqueness; see Section~\ref{sec-common}. It was also needed in
\cite{GorMar20-imm} to prove the following general inequalities extending
\eqref{cL-p}: for $p\in[1,\infty)$, $\alpha\ge-1/2$, and every
$n\in\mathbb{Z}_{+}$,
\begin{equation}\label{A-nn}
n^{(2\alpha+2)/p}\le A_{p\alpha}(n)\le
(n+\theta_{p\alpha})^{(2\alpha+2)/p},
\end{equation}
where
\[
A_{p\alpha}(n)=
\frac{\mathcal{C}_{p\infty,|\!\sin x|^{2\alpha+1}}(\mathcal{T}_{n})}
{\mathcal{C}_{p\infty,|x|^{2\alpha+1}}(\mathcal{E}_{1})}
=
\frac{\mathcal{C}_{p\infty,(1-x^{2})^{\alpha}}(\mathcal{P}_{n})}
{\mathcal{C}_{p\infty,x^{2\alpha+1};\mathbb{R}_{+}}
(\mathcal{E}_{1}^{\textup{even}})},\quad
\theta_{p\alpha}=\begin{cases}
\lceil\frac1p\rceil,&\alpha=-1/2,\\
2\lceil\frac{\alpha+3/2}{p}\rceil,&\alpha>-1/2.
\end{cases}
\]
These bounds imply Levin--Lubinsky type asymptotic formulas with an estimate of
the remainder; see~\cite{Gan19}.

The algebraic constant on $[-1,1]$ and the constant for even functions on
$\mathbb{R}_{+}$ that appear here follow directly from~\eqref{C-0}. In the
unweighted trigonometric case $\alpha=-1/2$ and the unweighted algebraic case
$\alpha=0$, the condition $p\ge1$ can be weakened to $p>0$; see
\cite{GorMar20-imm}. For other $\alpha$, this has been proved only when the
value at a fixed point is used explicitly in place of the $L^{\infty}$ norm.
In other words, extending~\eqref{A-nn} to $p<1$ requires an analogue of the
Arestov--Deikalova maximum-shift result. H\"older's inequality, which is used
for $p\ge1$, is no longer available. Since~\eqref{A-nn} is quite general, we
give its derivation in Section~\ref{sec-ll}.

Using the relation between multidimensional and one-dimensional weighted
Nikolskii constants for $\alpha=d/2-1$~\cite{GorMar20-imm}, we obtain the
following universal inequalities for $d\ge1$, $p\ge1$, and $n\ge0$:
\begin{equation}\label{n-C-l-n}
n^{d/p}\le
\frac{\mathcal{C}_{p\infty}(\Pi_{n}^{d})}{\mathcal{L}_{dp}}
\le(n+\theta_{p,d/2-1})^{d/p}.
\end{equation}
For $p=1$ these estimates can be continued using~\eqref{LA-d}. It would be
interesting to obtain similar bounds for multidimensional trigonometric
polynomials with spectrum in a ball.

The theory of Bernstein--Markov--Nikolskii constants continues to develop and
still contains many open questions, some of which are listed in the
conclusion.

\section{Basic methods for estimating Bernstein--Nikolskii constants}\label{sec-common}

We collect some general facts about estimates of Bernstein--Nikolskii (and
Markov--Nikolskii) constants that have been used, in one form or another, by
many authors; see, for example,~\cite{AreDei16,GorDob18,Gan19}. Not
surprisingly, many of these facts parallel classical questions in approximation
theory, including existence, uniqueness, characterization, and rates of
approximation; see, for example,~\cite{Tim60,Sha71,Kor76}. As above, let
\[
\mathcal{C}_{pq}^{D}(Y)=\sup_{f\in Y\cap L^{p}(M)}
\frac{\|Df\|_{q}}{\|f\|_{p}}.
\]
We restrict ourselves to $0<p\le q\le\infty$ (see the motivation in
Subsection~\ref{subsec-hist-nik}), assume that $Y$ is a subspace of $C(M)$,
and that $D$ is a linear differential operator with $DY\subset Y$; the
$L^{p}$ spaces may be weighted. The problem is to estimate, or whenever
possible compute, $\mathcal{C}_{pq}^{D}(Y)$ and at the same time verify that it
is finite. A simple example, $D=\partial$ and
$Y=\operatorname{span}\{\sin x^{2}\}\subset L^{\infty}(\mathbb{R})$, shows
that finiteness is not automatic.

\subsection[Reduction of pq inequalities to p-infinity and pp inequalities]{Reduction of $pq$ inequalities to $p\infty$ and $pp$ inequalities}\label{subsec-common-reduct}
Let $f\in Y$. We have
\[
\|Df\|_{q}\le \|Df\|_{\infty}^{1-p/q}\|Df\|_{p}^{p/q}.
\]
Hence
\[
\mathcal{C}_{pq}^{D}(Y)\le
(\mathcal{C}_{p\infty}^{D}(Y))^{1-p/q}
(\mathcal{C}_{pp}^{D}(Y))^{p/q},
\]
and, in particular, for $D=I$,
\begin{equation}\label{pq}
\mathcal{C}_{pq}(Y)\le(\mathcal{C}_{p\infty}(Y))^{1-p/q}.
\end{equation}
For example, if $Y=Y_{n}$, $\dim Y_{n}\asymp n^{d}$, and
$\mathcal{C}_{p\infty}^{D}(Y_{n})\le Cn^{r+d/p}$,
$\mathcal{C}_{pp}^{D}(Y_{n})\le Cn^{r}$, then we obtain an upper estimate of
the form~\eqref{C-asymp}:
\begin{equation}\label{C-n}
\mathcal{C}_{pq}^{D}(Y_{n})\le Cn^{r+d(1/p-1/q)}.
\end{equation}

Another way to derive a $pq$ Bernstein--Nikolskii inequality from a $qq$
Bernstein inequality and a $pq$ Nikolskii inequality is
\[
\|Df\|_{q}\le \mathcal{C}_{qq}^{D}(Y)\|f\|_{q}\le
\mathcal{C}_{qq}^{D}(Y)\mathcal{C}_{pq}(Y)\|f\|_{p},
\]
so by~\eqref{pq},
\[
\mathcal{C}_{pq}^{D}(Y)\le
\mathcal{C}_{qq}^{D}(Y)\mathcal{C}_{pq}(Y)\le
\mathcal{C}_{qq}^{D}(Y)(\mathcal{C}_{p\infty}(Y))^{1-p/q}.
\]
For the example above, this again gives~\eqref{C-n}. Unfortunately, this method
does not usually produce sharp constants, because the extremal functions in
the sharp Bernstein and Nikolskii inequalities are different.

\subsection{A multiplicative estimate for Nikolskii constants}
Let $D=I$ and suppose that the subspaces $Y=Y_{n}$, $n\ge0$, have the
multiplicative property: if $f\in Y_{n}$, then $f^{s}\in Y_{sn}$ for
$s\in\mathbb{Z}_{+}$. This is true, for example, for polynomials and entire
functions of exponential type. Then for every $f\in Y_{n}$,
\[
\|f\|_{\infty}=\|f^{s}\|_{\infty}^{1/s}\le
(\mathcal{C}_{p/s}(Y_{sn})\|f^{s}\|_{p/s})^{1/s}
=(\mathcal{C}_{p/s}(Y_{sn}))^{1/s}\|f\|_{p}.
\]
Therefore
\begin{equation}\label{mult}
\mathcal{C}_{p}(Y_{n})\le
(\mathcal{C}_{p/s}(Y_{sn}))^{1/s},\quad
s\in\mathbb{Z}_{+},\quad p>0.
\end{equation}

For homogeneous polynomials or spherical harmonics of order~$n$, for example,
this multiplicative property fails. This makes the estimation of
Bernstein--Nikolskii constants on these subspaces substantially more difficult.

It is also difficult to extend the multiplicative estimate to
Bernstein--Nikolskii constants with $D\ne I$, because the argument would require
an estimate of the form
$\|Df\|_{\infty}\le C\|Df^{s}\|_{\infty}^{1/s}$.

\subsection{Existence of an extremal function}\label{subsec-common-exist}
By linearity of $D$,
\begin{equation}\label{C-B}
\mathcal{C}_{pq}^{D}(Y)=\sup_{f\in B}\|Df\|_{q},\quad
B=\{f\in Y\colon \|f\|_{p}\le1\}.
\end{equation}

Suppose that $\mathcal{C}_{pq}^{D}(Y)<\infty$. Then the functional
$\|Df\|_{q}$ is continuous on the set (ball) $B$. If $\dim Y<\infty$, as in
the polynomial case, existence of an extremal function follows from
compactness of~$B$.

For $\dim Y=\infty$, we consider only the example of entire functions of
exponential type. Let $Y=\mathcal{E}_{\sigma}^{p}(\Omega)$, $p<\infty$. We use
the following compactness statement~\cite{Nik75}: if
$\{f_{n}\}_{n=1}^{\infty}\subset\mathcal{E}_{\sigma}(\Omega)$ is a sequence
of entire functions of exponential type uniformly bounded on
$\mathbb{R}^{d}$ by the same constant $C$, then it has a subsequence
$\{f_{n_{j}}\}$ converging uniformly on every compact subset of
$\mathbb{R}^{d}$ to a function $f\in\mathcal{E}_{\sigma}(\Omega)$, which is
bounded on $\mathbb{R}^{d}$ by the same constant $C$.

We now restrict to $q=\infty$ and omit $\infty$ from the subscript. Assume
that $D$ commutes with translations and is continuous with respect to locally
uniform convergence. Then~\eqref{C-B} gives
\begin{equation}\label{C-inf}
\mathcal{C}_{p}^{D}(Y)=\sup_{f\in B}\|Df\|_{\infty}.
\end{equation}
Let $\{f_{n}\}_{n=1}^{\infty}\subset B$ be a maximizing sequence with
$\|Df_{n}\|_{\infty}\to\mathcal{C}_{p}^{D}(Y)$. Choose points $x_{n}$ such
that
\[
|Df_{n}(x_{n})|\ge\|Df_{n}\|_{\infty}-\frac1n,
\]
and, using translation invariance of the space and norm, replace $f_{n}(x)$ by
$f_{n}(x+x_{n})$. Then
\[
|Df_{n}(0)|\to\mathcal{C}_{p}^{D}(Y).
\]
Since $\|f_{n}\|_{\infty}\le C\|f_{n}\|_{p}\le C$, there is a subsequence
$\{f_{n_{j}}\}$ converging as above to a function $f$. For every $R>0$,
\[
\|f\chi_{|x|\le R}\|_{p}
=\lim_{j\to\infty}\|f_{n_{j}}\chi_{|x|\le R}\|_{p}\le1,
\]
where $\chi_{A}$ is the characteristic function of $A$. Letting $R\to\infty$
gives $f\in B$. Moreover, locally uniform convergence and the assumptions on
$D$ imply
\[
Df_{n_{j}}(0)\to Df(0).
\]
Consequently,
\[
|Df(0)|=\mathcal{C}_{p}^{D}(Y),
\]
so $f$ is an extremal function.

\subsection{Shifting the maximum}\label{subsec-common-shift}
We now discuss when the $L^{\infty}$ norm in~\eqref{C-inf} can be replaced by
the value at a fixed point $o\in M$; see also~\cite{Gan19}. Introduce the
shifted Bernstein--Nikolskii problem
\[
\mathcal{A}_{p}^{D}(Y;o)=\sup_{f\in Y\cap L^{p}(M)}
\frac{|Df(o)|}{\|f\|_{p}}.
\]
If $Y=-Y$, as is the case for real linear subspaces, then in the real setting
the absolute value in $|Df(o)|$ can be omitted. In the complex setting, if $Y$
is invariant under multiplication by unimodular constants, the absolute value
can be replaced by $\operatorname{Re}Df(o)$.

For every subspace $Y_{0}\subseteq Y$,
$\mathcal{A}_{p}^{D}(Y_{0};o)\le\mathcal{C}_{p}^{D}(Y)$. In many cases,
however, there is a subspace $Y_{0}$ such that
$\mathcal{C}_{p}^{D}(Y)=\mathcal{A}_{p}^{D}(Y_{0};o)$. We give several
examples.

The simplest case is when $M$ is an abelian group under addition, $o=0$, and
$Y$, $D$, and the $p$-norm are invariant under translations $f(x+t)$, as in
the unweighted setting. Without assuming the existence of an extremal
function, let $\varepsilon>0$ and choose $f_{\varepsilon}\in B$ such that
$\|Df_{\varepsilon}\|_{\infty}>
\mathcal{C}_{p}^{D}(Y)-\varepsilon/2$. Choose $x_{\varepsilon}\in M$ with
\[
|Df_{\varepsilon}(x_{\varepsilon})|>
\|Df_{\varepsilon}\|_{\infty}-\varepsilon/2.
\]
Then
$|Df_{\varepsilon}(x_{\varepsilon})|>
\mathcal{C}_{p}^{D}(Y)-\varepsilon$. Therefore
$F(x)=f_{\varepsilon}(x+x_{\varepsilon})\in B$,
$|DF(0)|=|Df_{\varepsilon}(x_{\varepsilon})|$, and
\[
\mathcal{A}_{p}^{D}(Y;0)\le\mathcal{C}_{p}^{D}(Y)
<|DF(0)|+\varepsilon\le
\mathcal{A}_{p}^{D}(Y;0)+\varepsilon.
\]
Thus $\mathcal{C}_{p}^{D}(Y)=\mathcal{A}_{p}^{D}(Y;0)$. Note that we do not
assume here that the $L^{\infty}$ norm of an arbitrary function in $Y$ is
attained at some point. This is automatic when $M$ is compact. For entire
functions in $\mathcal{E}_{\sigma}^{p}(\Omega)$ with $p<\infty$, one may use
the fact that they tend to zero on $\mathbb{R}^{d}$ as $|x|\to\infty$. For
$1\le p<\infty$ this follows from~\cite{Nik75}; for $p<1$ we use the embedding
into $\mathcal{E}_{\sigma}^{1}(\Omega)$.

If $M$ is a homogeneous space with a transitive motion group $G$ and a
$G$-invariant measure, the same argument applies, provided that $Y$, $D$, and
the $p$-norm are invariant under $f(gx)$, $g\in G$. The point $o\in M$ can be
chosen arbitrarily and regarded as the origin on $M$. On the sphere
$\mathbb{S}^{d}$, for example, we may take the north pole
$o=e_{d+1}=(0,\ldots,0,1)$. To shift the maximum to $o$, set
$F(x)=f_{\varepsilon}(g_{\varepsilon}x)$, where the motion
$g_{\varepsilon}$ is chosen so that $g_{\varepsilon}o=x_{\varepsilon}$.

Suppose now that neither of these situations applies. For example, let
$M=\mathbb{T}$ and let $Y$ be the space of even trigonometric polynomials of
order at most $n$. Then $Y$ is not invariant under ordinary translations
$f(x+t)$. In such cases invariant generalized translation operators
$T^{t}f(x)$ are often useful. For their general theory see, for example,
\cite{Lev73}. That source, however, contains little about $L^{p}$ boundedness
and only a few examples related to singular Sturm--Liouville problems. Many
important examples of generalized translation operators with the properties
needed for sharp Nikolskii constants can be found in
\cite{AreDei15,AreDei16,AreBabDeiHor18,GorDob18,GorIva19,GorIvaTik19,Mar20};
these papers also contain further references.

We illustrate the idea for $p\ge1$ using
$Y=\mathcal{T}_{n}^{\text{even}}$. For even functions, a natural invariant
generalized translation is the symmetric operator
$T^{t}f(x)=\frac{f(x+t)+f(x-t)}{2}$. This choice is motivated, among other
things, by its action on the basis functions:
$T^{t}\cos kx=\cos kt\cos kx$. This is an instance of the general product
formula, which is often used to construct an integral representation of the
operator $T$. Let $D=\partial^{r}$. The operator $D$ is invariant under $T$
only for even $r$; for general $r$ it should be defined as a multiplier. If,
as above,
$\|Df_{\varepsilon}\|_{\infty}=|Df_{\varepsilon}(x_{\varepsilon})|$, set
$F(t)=T^{t}f_{\varepsilon}(x_{\varepsilon})$. Then
$DF(0)=Df_{\varepsilon}(x_{\varepsilon})$, $F\in Y$, and for $p\ge1$ the
triangle inequality gives
$\|F\|_{p}\le\|f_{\varepsilon}\|_{p}\le1$. Hence $F\in B$ and
$\mathcal{C}_{p}^{D}(Y)=\mathcal{A}_{p}^{D}(Y;0)$. The case $p<1$ is much
more difficult.

If instead $Y=\mathcal{T}_{n}$ is the full polynomial space, the same argument
gives, for even $r$,
\[
\mathcal{C}_{p}^{\partial^{r}}(\mathcal{T}_{n})=
\mathcal{C}_{p}^{\partial^{r}}(\mathcal{T}_{n}^{\text{even}})=
\mathcal{A}_{p}^{\partial^{r}}(\mathcal{T}_{n}^{\text{even}};0).
\]
For odd $r$, one may use odd polynomials and the nonsymmetric generalized
translation
$T^{t}f(x)=\frac{f(x+t)-f(x-t)}{2}$, which is invariant on that class. Similar
relations between Bernstein--Nikolskii constants on different classes hold in
more general settings. In particular, this method can reduce multidimensional
problems to one-dimensional weighted ones; see, for example,
\cite{Iva92,AreDei13}. At present, however, such reductions are available only
for $p\ge1$.

Consider one more example: a rank-one homogeneous space $M$ with metric
$\rho$, assumed compact for simplicity. Let $Y_{n}$ be the space of
polynomials of order at most $n$. As generalized translation, take the mean
operator
$T^{t}f(x)=\frac{1}{w(t)}\int_{\rho(x,\xi)=t}f(\xi)\,d\xi$, where
$w(t)=|\{\xi\colon\rho(x,\xi)=t\}|$ is the zonal weight, independent of $x$.
We may assume $t\in[0,\pi]$. For simplicity, let $D=I$, although powers of the
Laplace--Beltrami operator can also be used as in
Subsection~\ref{subsec-hist-nik}. The operator $T$ is nonsymmetric here. Let
$o\in M$ be the origin. We identify a zonal function
$f(x)=f(\rho(x,o))$ both with a function on $M$ and with the corresponding
one-dimensional function $f(t)$. Set $F(t)=T^{t}f(o)$ for $f\in Y_{n}$. Then
the zonal function $F(x)=F(\rho(x,o))$ belongs to $Y_{n}$, as follows from
harmonic analysis on $M$; see~\cite{Vil91}. The zonal polynomials form a
subspace $Y_{n0}\subset Y_{n}$, which can be identified with the even
trigonometric polynomials $\mathcal{T}_{n}^{\text{even}}$, or, after the
trigonometric change of variables, with the algebraic polynomials
$\mathcal{P}_{n}$.

We have the integral formula
\[
\int_{M}f(x)\,dx=a\int_{0}^{\pi}T^{t}f(o)w(t)\,dt,\quad
a=\frac{|M|}{\int_{0}^{\pi}w(t)\,dt},
\]
and H\"older's inequality gives, for $p\ge1$,
\begin{equation}\label{int-T}
\|F\|_{p}^{p}=a\int_{0}^{\pi}|F(t)|^{p}w(t)\,dt
=\|F\|_{p,aw}^{p}\le
a\int_{0}^{\pi}T^{t}(|f|^{p})(o)w(t)\,dt=\|f\|_{p}^{p},
\end{equation}
where the unweighted norm is taken in $L^{p}(M)$.

Now suppose, as above, that
$\|f_{\varepsilon}\|_{\infty}=|f_{\varepsilon}(x_{\varepsilon})|$, and set
$F(t)=T^{t}f_{\varepsilon}(x_{\varepsilon})$. Then
$F\in\mathcal{T}_{n}^{\text{even}}$ and, at the same time,
$F\in Y_{n0}$ is a zonal polynomial. We have
$F(0)=F(o)=f_{\varepsilon}(x_{\varepsilon})$. By~\eqref{int-T},
$\|F\|_{p}=\|F\|_{p,aw}\le\|f_{\varepsilon}\|_{p}\le1$. Thus, for
$p\ge1$,
\[
\mathcal{C}_{p}(Y_{n})=
\mathcal{C}_{p}(Y_{n0})=
\mathcal{C}_{p,aw}(\mathcal{T}_{n}^{\text{even}})=
\mathcal{A}_{p,aw}(\mathcal{T}_{n}^{\text{even}};0),
\]
which relates multidimensional constants to one-dimensional weighted ones. The
same method also applies to locally compact spaces $M$. For
$\mathbb{R}^{d}$, for example,
$w(t)=|\mathbb{S}^{d-1}|t^{d-1}$, $t\in\mathbb{R}_{+}$; the zonal functions are
radial functions $F(|x|)$, $Y_{n}$ is the space of entire functions of
spherical exponential type at most $n$, and $Y_{n0}$ can be identified with
even one-dimensional entire functions of exponential type at most $n$. For
harmonic analysis on rank-one hyperbolic spaces see~\cite{Hel64,Koo84}.

When $M$ carries a weighted $L^{p}$ space, the situation is more difficult.
One must establish the existence of invariant generalized translations with
the required $p$-norm estimates. The relevant tools use integral
representations of generalized translations together with invariance,
positivity, locality, and adjointness properties. Some progress has been made
in the papers cited above, which treat symmetric and nonsymmetric positive
Gegenbauer, Jacobi, Bessel, Gegenbauer--Dunkl, Bessel--Dunkl, and Dunkl
generalized translations.

\subsection{Reduction to the real case}\label{subsec-common-real}
Shifting the maximum makes it possible to reduce problems for complex
polynomials and entire functions of exponential type to real ones. Suppose,
for example, that
$\mathcal{C}_{p}^{D}(Y)=\mathcal{A}_{p}^{D}(Y;0)$, the subspace $Y$ is
invariant under
\[
Rf(z)=\frac{f(z)+\overline{f(\overline{z})}}{2},
\]
as is typical for analytic function spaces, and $D$ commutes with $R$. Choose
$f_{\varepsilon}\in B$ such that
$|Df_{\varepsilon}(0)|>\mathcal{C}_{p}^{D}(Y)-\varepsilon$. Without loss of
generality, assume that $Df_{\varepsilon}(0)\in\mathbb{R}_{+}$. Then
$f=Rf_{\varepsilon}$ is real-valued and belongs to $Y$. Moreover,
$Df(0)=Df_{\varepsilon}(0)$ and
$|f(x)|\le|f_{\varepsilon}(x)|$ for $x\in\mathbb{R}$, hence
$\|f\|_{p}\le\|f_{\varepsilon}\|_{p}$ for $p>0$. This proves the desired
reduction.

\subsection{Uniqueness of an extremal function}
Continue with the setting of Subsection~\ref{subsec-common-real}. Shifting the
maximum also proves uniqueness of an extremal function when $1<p<\infty$,
because $L^{p}$ is strictly convex. Suppose that there are two distinct
extremal functions $f_{1}$ and $f_{2}$ with
$Df_{1}(0)=Df_{2}(0)$ (the absolute value has been removed) and
$\|f_{1}\|_{p}=\|f_{2}\|_{p}$. Set
$f=\frac{f_{1}+f_{2}}{2}$. Then $f\in Y$, $Df(0)=Df_{1}(0)$, and strict
convexity gives $\|f\|_{p}<\|f_{1}\|_{p}$. Hence
$\frac{|Df(0)|}{\|f\|_{p}}>
\frac{|Df_{1}(0)|}{\|f_{1}\|_{p}}$, contradicting the extremality of $f_{1}$.

The cases $p=1,\infty$ are usually handled using a characterization of the
extremal functions. For $p=1$, for example, uniqueness can be proved from the
orthogonality relation~\eqref{h-f-0} in Subsection~\ref{subsec-common-orth};
see, for example,~\cite{DaiGorTik21}. One should also keep in mind that an
extremal function may be unique on a subset $Y_{0}\subset Y$ in the shifted
problem $\mathcal{A}_{p}^{D}(Y_{0};o)$ that yields
$\mathcal{C}_{p}^{D}(Y)$ (see Subsection~\ref{subsec-common-shift}), while its
uniqueness on the full class $Y$ may remain unclear.

\subsection{Basic estimates based on a reproducing kernel}\label{subsec-common-repr}
We present a general method for obtaining estimates of type \eqref{nw-nik} for
Nikolskii constants ($D=I$). It is based on estimating the reproducing kernel of
the subspace $Y$. A function $K(x,y)$ is such a kernel if it belongs to $Y$ in
both variables and
\begin{equation}\label{f-K}
f(x)=\int_{M}K(x,y)f(y)\,dy,\quad \forall\,f\in Y.
\end{equation}

For example, for $M=\mathbb{T}$ and $Y=\mathcal{T}_{n}$ we have
$f(x)=\int_{\mathbb{T}}D_{n}(x-y)f(y)\,dy$, where $D_{n}$ is the Dirichlet
kernel:
\[
D_{n}(x)=\frac{1}{2\pi}\sum_{|k|\le n}e^{ikx}=\frac{\sin((n+\frac{1}{2})x)}{2\pi \sin
(\frac{x}{2})}.
\]
If $M$ is a compact rank-one Riemannian manifold of dimension $d$ and
$Y=\Pi_{n}(M)$, then $K(x,y)$ is expressed in terms of the zonal Jacobi
polynomial $P_{n}^{(\alpha_{d},\beta_{d})}(\cos \rho(x,y))$. See, for example,
\cite{ConSlo99}, where the values of the parameters $\alpha_{d}$ and
$\beta_{d}$ are given (for the sphere, for instance,
$\alpha_{d}=d/2$, $\beta_{d}=d/2-1$). In the case $M=\mathbb{R}^{d}$ and
$Y=\mathcal{E}_{\sigma}^{2}(\Omega)$, we have
$K(x,y)=\widehat{\chi}_{\sigma\Omega}(x-y)$. In these examples the integral
operator in \eqref{f-K} is the orthogonal projector from $L^{2}(M)$ onto the
subspace $Y$.

Let $o\in M$ be the origin (as in Subsection~\ref{subsec-common-shift}) and put
$K(x)=K(o,x)$. In the examples above, the reproducing kernel has the property
that for every $x\in M$
\begin{equation}\label{K-K}
\|K(x,\cdot)\|_{2}^{2}=K(x,x)=\sup_{x,y\in M}|K(x,y)|,
\end{equation}
in particular, $\|K\|_{2}^{2}=\|K\|_{\infty}=K(o)$. The first equality in
\eqref{K-K} follows from \eqref{f-K} and the symmetry of the reproducing
kernel, while the second equality for homogeneous spaces $M$ follows from
harmonic analysis on~$M$ (see~\cite{Vil91}). In the general case it has to be
verified separately. We also note that for finite-dimensional invariant
subspaces $Y$ on a homogeneous space $M$ one has
$\|K\|_{2}^{2}=C\dim Y$; for example,
$\|D_{n}\|_{2}^{2}=\frac{\dim \mathcal{T}_{n}}{2\pi}$.

Assume that \eqref{K-K} holds. Then for $0<p\le 2$, the Cauchy--Schwarz
inequality gives
\[
\|f\|_{\infty}=\sup_{x}\,\biggl|\int_{M}K(x,\cdot)f\,dy\biggr|\le
\|f\|_{\infty}^{1-p/2} \sup_{x}\,\int_{M}|K(x,\cdot)|\,|f|^{p/2}\,dy\le
\|f\|_{\infty}^{1-p/2}\|K\|_{2}\|f\|_{p}^{p/2},
\]
whence $\|f\|_{\infty}\le (\|K\|_{2}^{2})^{1/p}\|f\|_{p}$ and
\begin{equation}\label{p-1-2}
\mathcal{C}_{p}(Y)\le (\|K\|_{2}^{2})^{1/p},\quad 0<p\le 2.
\end{equation}
This estimate is sharp for $p=2$, since by \eqref{K-K} the reproducing kernel
$K\in Y$ itself is an extremal function.

For $1\le p\le 2$, one can also apply H\"older's inequality directly and
estimate the Bernstein--Nikolskii constants for multiplier-type operators $D$.
We have
\begin{equation}\label{D-K}
\|Df\|_{\infty} =\sup_{x}\,\biggl|\int_{M}D_xK(x,\cdot)f\,dy\biggr| \le
\|D_xK(x,\cdot)\|_{p'}\|f\|_{p},
\end{equation}
where $\|D_xK(x,\cdot)\|_{p'}=\|DK\|_{p'}$ for homogeneous spaces $M$.
This norm is difficult to handle directly. Since $DK\in Y$, the reproducing
property of the kernel and the Cauchy--Schwarz inequality give
\[
\|DK\|_{\infty}\le \|K\|_{2}\|DK\|_{2}.
\]
Therefore,
\[
\|DK\|_{p'}\le
\|DK\|_{\infty}^{1-2/p'}\|DK\|_{2}^{2/p'}
\le
\|K\|_{2}^{2/p-1}\|DK\|_{2}.
\]
Combining this with \eqref{D-K}, we obtain the following analogue of
\eqref{p-1-2}:
\[
\mathcal{C}_{p}^{D}(Y)\le
\|K\|_{2}^{2/p-1}\|DK\|_{2},\quad 1\le p\le 2.
\]
This estimate is sharp for $p=2$. For $D=I$ it reduces to \eqref{p-1-2}.

Inequality \eqref{p-1-2} extends to all $p>0$ for subspaces $Y=Y_{n}$ with the
multiplicative property (all the examples above have this property). Let
$K=K_{n}$ and $I_{n}=\|K_{n}\|_{2}^{2}$. We use \eqref{mult} and set
$s=p_{0}=\lceil \frac{p}{2}\rceil$. Then $p/s\le 2$, and hence
\begin{equation}\label{C-p-I}
\mathcal{C}_{p}(Y_{n})\le (\mathcal{C}_{p/s}(Y_{sn}))^{1/s}\le
I_{p_{0}n}^{1/p}.
\end{equation}
Using \eqref{pq}, we also obtain
\[
\mathcal{C}_{pq}(Y_{n})\le (\mathcal{C}_{p}(Y_{n}))^{1-p/q}\le
I_{p_{0}n}^{1/p-1/q}.
\]
For example, for the Dirichlet kernel we have
$I_{n}=\frac{2n+1}{2\pi}$, and therefore
$\mathcal{C}_{pq}(\mathcal{T}_{n})\le
(\frac{2p_{0}n+1}{2\pi})^{1/p-1/q}$, which agrees with \eqref{C-T-d} for
$d=1$. In all the examples above, $I_{n}\asymp n^{d}$. For finite-dimensional
$Y_{n}$ this also follows from harmonic analysis. Thus
$\mathcal{C}_{pq}(Y_{n})\le C(p_{0}n)^{d(1/p-1/q)}$, where the constant $C$
can be estimated effectively. This yields basic inequalities of type
\eqref{nw-nik}.

To prove that estimate \eqref{C-p-I} has the correct order, one usually uses
Fej\'er--Jackson type kernels $F_{n}=K_{n/s}^{s}$, where
$s=\lceil \frac{2}{p}\rceil$ (for simplicity, assume that the integer $n$ is
divisible by~$s$). Then $F_{n}\in Y_{n}$,
$\|F_{n}\|_{\infty}=\|K_{n/s}\|_{\infty}^{s}$, and
\[
\|F_{n}\|_{p}=\biggl(\int_{M}|K_{n/s}|^{2}|K_{n/s}|^{sp-2}\,dx\biggr)^{1/p}\le
\|K_{n/s}\|_{\infty}^{s-2/p}\|K_{n/s}\|_{2}^{2/p}.
\]
Since $\|K_{n/s}\|_{\infty}=\|K_{n/s}\|_{2}^{2}$, it follows that
\[
\frac{\|F_{n}\|_{\infty}}{\|F_{n}\|_{p}}\ge
\|K_{n/s}\|_{\infty}^{2/p}\|K_{n/s}\|_{2}^{-2/p}
=(\|K_{n/s}\|_{2}^{2})^{1/p}.
\]

Thus,
\[
I_{n/s}^{1/p}\le \mathcal{C}_{p}(Y_{n})\le I_{p_{0}n}^{1/p},\quad s=\lceil
\tfrac{2}{p}\rceil,\quad p_{0}=\lceil \tfrac{p}{2}\rceil.
\]
For $p=2$ these estimates coincide. For $p=1$ we have
\begin{equation}\label{C-p-1}
I_{n/2}\le \mathcal{C}_{1}(Y_{n})\le I_{n}.
\end{equation}
If $I_{n}\asymp n^{d}$ as $n\to \infty$, then
\[
\mathcal{C}_{p}(Y_{n})\asymp n^{d/p}.
\]

The same argument gives equivalent estimates for the shifted problem
$\mathcal{A}_{p}(Y;o)$. In this case properties of type \eqref{K-K} are needed
only at $x=o$.

We have considered homogeneous spaces $M$, which can be realized as quotient
spaces $M=G/H$, where $H$ is the stabilizer of the point $o$. In
\cite{NurRuzTik15} these estimates were extended to compact Lie groups $G$.
However, $G$ itself can be realized as the homogeneous space $(G\times G)/G$,
so the above scheme applies directly to groups as well.

The results also extend to weighted settings, where harmonic analysis is
related to generalized translation operators.

\subsection{The Nikolskii constant for nonnegative functions}\label{subsec-common-pos}
In many cases, for homogeneous spaces $M$ with a real zonal reproducing kernel
$K_{n}$, the lower estimate in \eqref{C-p-1} has an extremal interpretation.
Let $Y_{n0}$ be the subspace of zonal functions in $Y_{n}$. Consider the
shifted Nikolskii problem for nonnegative functions: find
\[
\mathcal{A}_{1}^{+}(Y_{n0};0)=\sup_{0\le f\in Y_{n0}\cap
L^{1}(M)}\frac{f(0)}{\|f\|_{1}}.
\]

In \cite{DaiGorTik20-anal}, Gaussian quadrature formulas on the class $Y_{n0}$
were used to solve this problem. For simplicity, assume that $n$ is even when
$\dim Y_{n}<\infty$. In the cases considered there, a Markov--Gauss quadrature
formula of the form
\begin{equation}\label{quadr-f}
\int_{M}f(\rho(x,o))\,dx=\gamma_{0}f(0)+\sum_{k=1}^{\infty}\gamma_{k}f(t_{k})
\end{equation}
holds, where the weights $\gamma_{k}>0$, the points $t_{k}$ are the zeros of
the kernel $K_{n/2}$, and the sum is finite when $\dim Y_{n}<\infty$. This
formula is exact for the Fej\'er kernel $F_{n}=K_{n/2}^{2}$, so
$\|K_{n/2}^{2}\|_{1}=\gamma_{0}K_{n/2}^{2}(0)$, which is equivalent to
$\|K_{n/2}\|_{2}^{2}=\gamma_{0}\|K_{n/2}\|_{\infty}^{2}$, and hence
$\gamma_{0}=\|K_{n/2}\|_{2}^{-2}$. For an arbitrary nonnegative function
$f\in Y_{n0}$, dropping the nonnegative sum in \eqref{quadr-f} gives
$\|f\|_{1}\ge \gamma_{0}f(0)$, and therefore
$\frac{f(0)}{\|f\|_{1}}\le \gamma_{0}^{-1}=\|K_{n/2}\|_{2}^{2}$. This estimate
is sharp for $F_{n}$. Thus,
\[
\mathcal{A}_{1}^{+}(Y_{n0};0)=I_{n/2}.
\]

The required quadrature formulas for zonal functions are known for the sphere
$\mathbb{S}^{d}$ and projective spaces (where they coincide with the classical
weighted Gauss--Markov quadrature formulas), for Euclidean space
$\mathbb{R}^{d}$, and for hyperbolic spaces (see
\cite{GorIva15,GorIva18,GorIva19-sb} and the references therein).

\subsection{Orthogonality relations}\label{subsec-common-orth}
Continuing the preceding example, for simplicity consider the problem
\[
\mathcal{A}_{p}(Y_{n0};0)=\sup_{f\in Y_{n0}\cap
L^{p}(M)}\frac{|f(0)|}{\|f\|_{p}},
\]
where $Y_{n0}$ is a closed subspace of zonal functions and $1\le p<\infty$.
Then we have the convex problem
\begin{equation}\label{A-conv}
(\mathcal{A}_{p}(Y_{n0};0))^{-1}=\inf\,\{\|f\|_{p}\colon f\in Y_{n0},\
f(0)=1\}.
\end{equation}
Assume that an extremal function $f_{*}\in Y_{n0}$ exists. Then
\eqref{A-conv} can be rewritten as
\[
0=\operatorname*{argmin}\,\{\|f_{*}+h\|_{p}\colon h\in H\},
\]
where $H\subset Y_{n0}$ is the closed subspace defined by $h(0)=0$. A necessary
and sufficient condition for the zero element to be a best approximation to
$f_{*}$ from a closed subspace is that
$|f_{*}|^{p-1}\operatorname{sign}f_{*}$ be orthogonal to $H$ (see
\cite{Sha71}). For $p=1$, one also has to require that $f_{*}$ not vanish on a
set of positive measure (which holds for analytic functions). Thus we obtain
the orthogonality relation
\begin{equation}\label{h-f-0}
\int_{M}h|f_{*}|^{p-1}\operatorname{sign}f_{*}\,dx=0,\quad \forall\,h\in H,
\end{equation}
which characterizes the extremal function $f_{*}$. For an arbitrary
$f\in Y_{n0}$, taking $h=f-f(0)f_{*}$ also gives the equivalent
characterization
\begin{equation}\label{f0-p-1}
f(0)=\|f_{*}\|_{p}^{-p}\int_{M}f|f_{*}|^{p-1}\operatorname{sign}f_{*}\,dx.
\end{equation}
Relations \eqref{h-f-0} and \eqref{f0-p-1} can also be written for weighted
spaces $L_{w}^{p}$ with a zonal weight $w$ on $[0,\pi]$ when
$\dim Y_{n0}<\infty$, and on $\mathbb{R}_{+}$ otherwise.

Orthogonality relations can be used to prove uniqueness of extremal functions
in the non-strictly convex case $p=1$; see, for example,
\cite{AmiZie76,AreDei15} for algebraic polynomials. In \cite{DaiGorTik21}, this
approach was used to prove uniqueness of the radial extremal function in the
Nikolskii problem $\mathcal{L}_{d1}$.

For general Bernstein--Nikolskii constants the argument is similar; see, for
example, \cite{GorMar20-cheb}.

\subsection{Connection with the Chebyshev problem}
We continue the example from Subsection~\ref{subsec-common-orth} in the case
$\dim Y_{n0}<\infty$ and with a zonal weight $w$. We may then take
$Y_{n0}=\mathcal{T}_{n}^{\text{even}}$. Introduce the weight
$W(t)=(1-\cos t)w(t)$ and consider the trigonometric Chebyshev problem in
$L_{W}^{p}([0,\pi])$:
\[
T_{*}=\operatorname*{argmin}\,\{\|T\|_{p,W}\colon T(t)=\cos nt+S(t),\ S\in
Y_{n-1,0}\}.
\]
As above, the extremal polynomial is characterized by the orthogonality
relation
\begin{equation}\label{S-T-0}
\int_{0}^{\pi}S|T_{*}|^{p-1}\operatorname{sign}T_{*}\,W\,dt=0,\quad
\forall\,S\in Y_{n-1,0}.
\end{equation}
But $SW=\widetilde{S}w$, where
$\widetilde{S}(t)=S(t)(1-\cos t)\in Y_{n0}$ and $\widetilde{S}(0)=0$.
Therefore \eqref{S-T-0} is equivalent to
$\int_{0}^{\pi}\widetilde{S}|T_{*}|^{p-1}\operatorname{sign}T_{*}\,w\,dt=0$
for every $\widetilde{S}\in H$. Together with \eqref{h-f-0}, this shows that
the polynomials $f_{*}$ in the Nikolskii problem and $T_{*}$ in the Chebyshev
problem satisfy the same orthogonality relation, which uniquely characterizes
the extremal polynomial. Hence they coincide up to a constant factor. A
trigonometric substitution reduces the Chebyshev problem to an algebraic one.
Further properties of the trigonometric and algebraic polynomials, including
information about their zeros, can be derived as in
\cite{AmiZie76,AreDei13}.

\subsection{Dual problems}\label{subsec-common-dual}
Consider the shifted convex Nikolskii problem \eqref{A-conv}. Duality is widely
used in convex analysis and its applications to approximation theory. The dual
problem can be obtained either from general results, such as those presented
in \cite{Sha71}, or directly from the identity
$f(0)=\int_{M}K_{n}(x)f(x)\,dx$ and the characterization of the extremal
function in \eqref{f0-p-1}. Let $Y_{n0}^{\bot}$ denote the subspace of
$L^{p'}(M)$ orthogonal to $Y_{n0}\cap L^{p}(M)$. Then, for every
$g\in Y_{n0}^{\bot}$, H\"older's inequality gives
\[
|f(0)|=\biggl|\int_{M}(K_{n}(x)-g(x))f(x)\,dx\biggr|\le
\|K_{n}-g\|_{p'}\|f\|_{p}.
\]
Hence
$\mathcal{A}_{p}(Y_{n0};0)\le \inf_{g\in
Y_{n0}^{\bot}}\|K_{n}-g\|_{p'}$. To show that this estimate is sharp, take
$g_{*}=K_{n}-c|f_{*}|^{p-1}\operatorname{sign}f_{*}$, where
$c=\|f_{*}\|_{p}^{-p}$. Then
\[
\|K_{n}-g_{*}\|_{p'}=c\||f_{*}|^{p-1}\|_{p'}=\|f_{*}\|_{p}^{-1}=\mathcal{A}_{p}(Y_{n0};0),
\]
and for every $f\in Y_{n0}$, by \eqref{f0-p-1},
\[
\int_{M}fg_{*}\,dx=
\int_{M}fK_{n}\,dx-c\int_{M}f|f_{*}|^{p-1}\operatorname{sign}f_{*}\,dx=f(0)-f(0)=0.
\]

Thus,
\[
\mathcal{A}_{p}(Y_{n0};0)=\inf_{g\in
Y_{n0}^{\bot}}\|K_{n}-g\|_{p'}=\|K_{n}-g_{*}\|_{p'}.
\]
In other words, the dual Nikolskii problem is equivalent to the best
approximation of the reproducing kernel by the orthogonal complement of the
subspace.

As an application of the dual problem, consider the constant
$L=\mathcal{C}_{1}(\mathcal{E}_{1}^{1})$ (see
Section~\ref{subsec-hist-now}). For $L$ one can use the shifted problem, and
hence in $L^{\infty}(\mathbb{R})$
\[
L=\inf_{g\in (\mathcal{E}_{1}^{1})^{\bot}}\Bigl\|\frac{\sin t}{\pi
t}-g(t)\Bigr\|_{\infty}.
\]
This problem is difficult, so we consider its restriction to the subspace of
cosine series in $(\mathcal{E}_{1}^{1})^{\bot}$ (orthogonality follows, for
example, from the Paley--Wiener theorem):
\begin{equation}\label{L-I}
L\le \inf_{a_{1},a_{2},\ldots\in \mathbb{R}}\,\biggl\|\frac{\sin t}{\pi
t}-\sum_{k=1}^{\infty}a_{k}\cos kt\biggr\|_{\infty}=I.
\end{equation}
The following solution of this problem was proposed in \cite{Gor05}. A cosine
series is uniquely determined on $[0,\pi]$, so choose $a_{k}^{*}$ so that
\[
h_{*}(t)=\frac{\sin t}{\pi t}-\sum_{k=1}^{\infty}a_{k}^{*}\cos kt\equiv c,
\quad t\in [0,\pi].
\]
Integrating over $[0,\pi]$, we find
$c=\frac{\operatorname{Si}\pi}{\pi^{2}}$. Thus $h_{*}(t)=c$ on
$[0,\pi]$, while outside this interval one verifies that
$|h_{*}(t)|\le c$ for $t\in [0,\infty)$ (this was done in
\cite{AndKonPop96,AndKonPop00}). If there were another function
$h(t)=\frac{\sin t}{\pi t}-\sum_{k=1}^{\infty}a_{k}\cos kt$ such that
$\|h\|_{\infty}<\|h_{*}\|_{\infty}=c$, then
$h_{*}(t)-h(t)=\sum_{k=1}^{\infty}(a_{k}-a_{k}^{*})\cos kt>0$ for
$t\in [0,\pi]$, which is impossible because the integral of the difference is
zero. Hence
$I=\frac{\operatorname{Si}\pi}{\pi^{2}}=0.18764\ldots.$

In \cite{DaiGorTik21}, a similar method, in a more complicated setting, was
used to solve the $d$-dimensional version of this problem, which reduces to the
extremal problem in $L^{\infty}(\mathbb{R}_{+})$
\begin{equation}\label{j-j}
\inf_{a_{1},a_{2},\ldots\in \mathbb{R}}\,
\biggl\|j_{\alpha+1}(t)-\sum_{k=1}^{\infty}a_{k}j_{\alpha}(\tfrac{q_{\alpha+1,k}}{q_{\alpha+1,1}}\,t)\biggr\|_{\infty},
\end{equation}
where $0<q_{\alpha 1}<q_{\alpha 2}<\ldots$ are the zeros of the normalized
Bessel function $j_{\alpha}$ and $\alpha=d/2-1$. For example, when $d=1$ we
have $j_{-1/2}(t)=\cos t$, $j_{1/2}(t)=\frac{\sin t}{t}$,
$q_{1/2,k}=\pi k$, and we recover problem \eqref{L-I} solved above. In
\cite{DaiGorTik21}, problem \eqref{j-j} was solved for all
$\alpha\ge -0.272$. Thus all dimensions $d\ge 2$ are covered, and this gives
the upper estimate in \eqref{LA-d} from Subsection~\ref{subsec-hist-now} (see
\cite{DaiGorTik21} for details).

\section{The Levin--Lubinsky estimate with explicit bounds}\label{sec-ll}

In this section we prove the explicit universal bounds \eqref{A-nn}, which
imply Levin--Lubinsky type asymptotic formulas for the sphere
$\mathbb{S}^{d}$, the torus $\mathbb{T}$ with the periodic Gegenbauer weight
$|\!\sin x|^{2\alpha+1}$, and the interval $[-1,1]$ with the algebraic
Gegenbauer weight $(1-x^{2})^{\alpha}$ (see
\cite{DaiGorTik20-anal,Gan19,GorMar20-imm}).

Recall (see \eqref{A-nn}) that these bounds have the form
\[
n^{(2\alpha+2)/p}\le A_{p\alpha}(n)\le (n+\theta_{p\alpha})^{(2\alpha+2)/p},\quad
n\in \mathbb{Z}_{+},\quad 1\le p<\infty,\quad \alpha\ge -1/2,
\]
where
\[
A_{p\alpha}(n)=\frac{\mathcal{C}_{p\infty,|\!\sin
x|^{2\alpha+1}}(\mathcal{T}_{n})}
{\mathcal{C}_{p\infty,|x|^{2\alpha+1}}(\mathcal{E}_{1})}=
\frac{\mathcal{C}_{p\infty,(1-x^{2})^{\alpha}}(\mathcal{P}_{n})}
{\mathcal{C}_{p\infty,x^{2\alpha+1};\mathbb{R}_{+}}(\mathcal{E}_{1}^{\textup{even}})},\quad
\theta_{p\alpha}=\begin{cases} \lceil \frac{1}{p}\rceil,&\alpha=-1/2,\\ 2\lceil
\frac{\alpha+3/2}{p}\rceil,&\alpha>-1/2.
\end{cases}
\]

The starting point is the use of the shifted Nikolskii problems
\[
A_{p\alpha}(n;0)=\frac{\mathcal{A}_{p,|\!\sin
x|^{2\alpha+1}}(\mathcal{T}_{n};0)}
{\mathcal{A}_{p,|x|^{2\alpha+1}}(\mathcal{E}_{1};0)}.
\]
Then, by the scheme described in Subsection~\ref{subsec-common-shift} (see
\cite{GorMar20-imm} for details), one proves that
$A_{p\alpha}(n)=A_{p\alpha}(n;0)$ for $p\ge 1$. Thus it is enough to estimate
$A_{p\alpha}(n;0)$. We do this for all $0<p<\infty$. For brevity, put
$\nu=2\alpha+1\ge 0$.

\subsection{Lower estimate}
For $n=0$ the lower estimate is obvious, so let $n\ge 1$ and let $f\ne 0$ be
an arbitrary function in $\mathcal{E}_{1}\cap
L_{|x|^{\nu}}^{p}(\mathbb{R})$. The power weight $|x|^{\nu}$ on the real line
is a special case of a Dunkl weight, and therefore the weighted Nikolskii
inequality \cite{GorIvaTik19} implies that $f$ is bounded on the real line.

Consider the Fej\'er kernel
\[
\varphi(x)=\Bigl(\frac{\sin(x/2)}{x/2}\Bigr)^{2},\quad
\widehat{\varphi}(y)=(1-|y|)_{+},\quad x,y\in \mathbb{R}.
\]
We have the identity
$\sum_{k\in \mathbb{Z}}\varphi(x+2\pi k)=1$ for all $x\in\mathbb{R}$. For
example, it follows from the Poisson summation formula:
\[
\sum_{k\in \mathbb{Z}}\varphi(x+2\pi k)= \sum_{k\in
\mathbb{Z}}\widehat{\varphi}(k)e^{ikx}= \widehat{\varphi}(0)=\sum_{k\in
\mathbb{Z}}\varphi(2\pi k)=1.
\]

Set $g(x)=f(nx)\varphi(x)$. Then $g\in \mathcal{E}_{n+1}$ and
$g(x)=O((1+x^{2})^{-1})$ on the real line. Hence the Fourier transform
$\widehat{g}$ is continuous on $\mathbb{R}$ and, by the Paley--Wiener theorem,
$\widehat{g}(y)=0$ for $|y|\ge n+1$. Therefore the periodization of $g$
satisfies
\[
\sum_{k\in \mathbb{Z}}g(x+2\pi k)=\sum_{|k|\le n}\widehat{g}(k)e^{ikx}=T(x)\in
\mathcal{T}_{n}.
\]

We estimate from above the weighted norm
\[
A=\int_{-\pi}^{\pi}|T(x)|^{p}|\!\sin x|^{\nu}\,dx= \int_{-\pi}^{\pi}
\biggl|\sum_{k\in \mathbb{Z}}g(x+2\pi k)\biggr|^{p}|\!\sin x|^{\nu}\,dx
\]
of the polynomial $T$, where
\[
\biggl|\sum_{k\in \mathbb{Z}}g(x+2\pi k)\biggr|^{p}\le \biggl(\,\sum_{k\in
\mathbb{Z}}|f(n(x+2\pi k))|\varphi(x+2\pi k)\biggr)^{p}=B.
\]

We claim that
\[
B\le \sum_{k\in \mathbb{Z}}|f(n(x+2\pi k))|^{p}.
\]
For $p<1$,
\[
B\le \sum_{k\in \mathbb{Z}}|f(n(x+2\pi k))|^{p}\varphi^{p}(x+2\pi k)\le \sum_{k\in
\mathbb{Z}}|f(n(x+2\pi k))|^{p}.
\]
For $p\ge 1$, H\"older's inequality gives
\[
B\le \sum_{k\in \mathbb{Z}}|f(n(x+2\pi k))|^{p}\biggl(\,\sum_{k\in
\mathbb{Z}}\varphi^{p'}(x+2\pi k)\biggr)^{p/p'}\le \sum_{k\in
\mathbb{Z}}|f(n(x+2\pi k))|^{p},
\]
where we used
\[
\biggl(\,\sum_{k\in \mathbb{Z}}\varphi^{q}(x+2\pi k)\biggr)^{1/q}\le \sum_{k\in
\mathbb{Z}}\varphi(x+2\pi k)=1,\quad q\ge 1.
\]

Thus,
\begin{align*}
A&\le \int_{-\pi}^{\pi}\sum_{k\in \mathbb{Z}}|f(n(x+2\pi k))|^{p}|\!\sin x|^{\nu}\,dx=
\frac{1}{n}\sum_{k\in \mathbb{Z}}\int_{2\pi n(k-1/2)}^{2\pi
n(k+1/2)}|f(x)|^{p}|\!\sin \tfrac{x}{n}|^{\nu}\,dx\\
&=\frac{1}{n}\int_{\mathbb{R}}|f(x)|^{p}|\!\sin \tfrac{x}{n}|^{\nu}\,dx\le
\frac{1}{n^{\nu+1}}\int_{\mathbb{R}}|f(x)|^{p}|x|^{\nu}\,dx.
\end{align*}
Together with
\[
T(0)=\sum_{k\in \mathbb{Z}}f(2\pi kn)\varphi(2\pi k)=f(0),
\]
this gives
\[
n^{(\nu+1)/p}\,\frac{f(0)}{\|f\|_{p,|x|^{\nu}}}\le \frac{T(0)}{\|T\|_{p,|\!\sin
x|^{\nu}}}\le \mathcal{A}_{p,|\!\sin x|^{\nu}}(\mathcal{T}_{n};0).
\]
Since $f$ is arbitrary, this proves the required lower estimate.

\subsection{Upper estimate}
Let $n\ge 0$ and let $T\in \mathcal{T}_{n}\setminus \{0\}$ be an arbitrary
trigonometric polynomial. Clearly, $T\in \mathcal{E}_{n}$ as well.

First suppose that $\nu>0$. Take another Fej\'er kernel
$\psi(x)=(\frac{\sin x}{x})^{2}$ and let $s\ge 1$ be an integer to be chosen
later. Set $g(x)=\psi^{s}(x)T(x)$. Then $g\in \mathcal{E}_{n+2s}$.

We estimate the weighted norm of $g$ from above. On the real line,
\[
|g(x)|^{p}|x|^{\nu}= |\psi^{s}(x)T(x)|^{p}|x|^{\nu}= |T(x)|^{p}|\!\sin x|^{\nu}
\Bigl|\frac{\sin x}{x}\Bigr|^{2sp-\nu}=|T(x)|^{p}|\!\sin
x|^{\nu}\psi^{sp-\nu/2}(x),
\]
and therefore
\begin{align*}
\int_{\mathbb{R}}|g(x)|^{p}|x|^{\nu}\,dx&= \sum_{k\in
\mathbb{Z}}\int_{2\pi(k-1/2)}^{2\pi(k+1/2)}|T(x)|^{p}|\!\sin
x|^{\nu}\psi^{sp-\nu/2}(x)\,dx\\
&=\sum_{k\in \mathbb{Z}}\int_{-\pi}^{\pi}|T(x+2\pi k)|^{p}|\!\sin(x+2\pi
k)|^{\nu}\psi^{sp-\nu/2}(x+2\pi k)\,dx\\
&=\int_{-\pi}^{\pi}|T(x)|^{p}|\!\sin x|^{\nu}\sum_{k\in \mathbb{Z}}\psi^{sp-\nu/2}(x+2\pi k)\,dx.
\end{align*}
Now choose the smallest $s$ such that $sp-\nu/2\ge 1$. Then
$s=\lceil\frac{1+\nu/2}{p}\rceil$ and
\[
\sum_{k\in \mathbb{Z}}\psi^{sp-\nu/2}(x+2\pi k)\le \biggl(\,\sum_{k\in
\mathbb{Z}}\psi(x+2\pi k)\biggr)^{sp-\nu/2}.
\]
Since $\psi$ is positive definite, its periodization is also positive
definite, and hence for all $x\in \mathbb{R}$
\[
\sum_{k\in \mathbb{Z}}\psi(x+2\pi k)\le \sum_{k\in \mathbb{Z}}\psi(2\pi
k)=\psi(0)=1.
\]

Thus,
\[
\int_{\mathbb{R}}|g(x)|^{p}|x|^{\nu}\,dx\le \int_{-\pi}^{\pi}|T(x)|^{p}|\!\sin
x|^{\nu}\,dx.
\]
Make the change of variables $g(x)=f((n+2s)x)$. Then
$f\in \mathcal{E}_{1}$,
\[
\int_{\mathbb{R}}|f(x)|^{p}|x|^{\nu}\,dx=
(n+2s)^{\nu+1}\int_{\mathbb{R}}|g(x)|^{p}|x|^{\nu}\,dx\le
(n+2s)^{\nu+1}\int_{-\pi}^{\pi}|T(x)|^{p}|\!\sin x|^{\nu}\,dx,
\]
and $f(0)=g(0)=T(0)$. Consequently,
\[
\frac{T(0)}{\|T\|_{p,|\!\sin x|^{\nu}}}\le
(n+2s)^{(\nu+1)/p}\,\frac{f(0)}{\|f\|_{p,|x|^{\nu}}}\le
(n+2s)^{(\nu+1)/p}\mathcal{A}_{p,|x|^{\nu}}(\mathcal{E}_{1};0).
\]
Since $T$ is arbitrary, this proves the required upper estimate.

If $\nu=0$, one can use the Fej\'er kernel
$(\frac{\sin(x/2)}{x/2})^{2}$. Then $n+2s$ is replaced by $n+s$, where
$s=\lceil \frac{1}{p}\rceil$.

\section{Applications of Bernstein--Nikolskii constants}\label{sec-appl}

Here we give several applications of Bernstein--Nikolskii constants in
different areas of mathematics. In particular, we formulate some results in
which these constants appear.

\subsection{Approximation theory}\label{subsec-appl-approx}
Order-sharp Bernstein--Nikolskii inequalities are widely used in approximation
theory, for example, in proofs of inverse theorems. Many results of this kind
are known. As an example, we refer to the recent papers
\cite{GorIvaTik19,GorIvaTik20-jat}, which study best approximation by entire
functions of spherical exponential type in $L_{v}^{p}(\mathbb{R}^{d})$ with a
power Dunkl weight $v$. In particular, the following versions of classical
inverse theorems of approximation theory were proved. Let
$1\le p\le \infty$, $r>0$, and $n\in \mathbb{N}$. Then for every
$f\in L_{v}^{p}(\mathbb{R}^{d})$,
\begin{equation}\label{omega-E}
\omega_{r}\Bigl(f,\frac{1}{n}\Bigl)_{p}\le C\frac{1}{n^{r}}\sum_{j=0}^{n}
(j+1)^{r-1}E_{j}(f)_{p},
\end{equation}
where $\omega_{r}(\cdot)_{p}$ is the modulus of smoothness of $f$ and
$E_{j}(\cdot)_{p}$ is its best approximation of order~$j$ by entire functions
of spherical exponential type at most $j$. For $1<p<\infty$ and
$q=\min\{p,2\}$,
\[
\omega_{r}\Bigl(f,\frac{1}{n}\Bigl)_{p}\le
C\,\frac{1}{n^{r}}\biggl(\sum_{j=1}^{n}j^{qr-1}
E_{j}^q(f)_{p}\biggr)^{1/q}+\frac{\|f\|_{p}}{n^{r}},
\]
which refines \eqref{omega-E}. The proofs use order-sharp Bernstein and
Nikolskii inequalities established in those papers in
$L_{v}^{p}(\mathbb{R}^{d})$ for entire functions of spherical exponential type
and fractional powers of the Dunkl Laplacian.

We also note that \cite{GorIvaTik20-jat} uses weighted Fourier--Dunkl
inequalities. These inequalities also estimate weighted $pq$-norms of
operators, which makes them related to Bernstein--Nikolskii inequalities. A
specific feature of the latter theory, however, is that one usually studies
functions whose Fourier transform has bounded support.

\subsection[A property of the Lp-approximant for p at least 1]{A property of the $L^{p}$-approximant for $p\ge 1$}\label{subsec-appl-jac}
In his pioneering paper \cite{Jac33}, D.~Jackson considered the following
problem for the torus and the interval: estimate, in the uniform norm, the
deviation of a continuous function from its polynomial of best approximation
(an approximant) taken in the $L^{p}$ metric, $1\le p<\infty$. The motivation
is as follows. Unlike in the uniform metric, an approximant in $L^{2}$ can be
constructed effectively, for example by the least-squares method. If it is
then used in the uniform metric, its error has to be estimated in that metric.
The following argument gives such an estimate. Let $M$ be a compact manifold
with normalized measure ($|M|=1$), and let $Y\subset C(M)$ be a polynomial
subspace. Denote by
\[
A_{p}f=\operatorname*{argmin}_{u\in Y\cap L^{p}(M)}\|f-u\|_{p}
\]
one of the approximants to $f$ in the $L^{p}$ metric (unique for
$1<p<\infty$). Then for $f\in C(M)$,
\begin{equation}\label{jac-res}
\|f-A_{p}f\|_{\infty}\le (2\mathcal{C}_{p}(Y)+1)\|f-A_{\infty}f\|_{\infty}.
\end{equation}

Indeed,
\[
\|f-A_{p}f\|_{\infty}\le
\|f-A_{\infty}f\|_{\infty}+\|A_{\infty}f-A_{p}f\|_{\infty}\le
\|f-A_{\infty}f\|_{\infty}+\mathcal{C}_{p}(Y)\|A_{\infty}f-A_{p}f\|_{p},
\]
where
\[
\|A_{\infty}f-A_{p}f\|_{p}\le \|f-A_{\infty}f\|_{p}+\|f-A_{p}f\|_{p}\le
2\|f-A_{\infty}f\|_{p}\le 2\|f-A_{\infty}f\|_{\infty}.
\]

If $\dim Y\asymp n^{d}$ and $\mathcal{C}_{p}(Y)\asymp n^{d/p}$ as
$n\to\infty$, then
$\|f-A_{p}f\|_{\infty}\le Cn^{d/p}\|f-A_{\infty}f\|_{\infty}$. Hence, if
$\|f-A_{\infty}f\|_{\infty}=o(n^{-d/p})$, the approximant $A_{p}f$ converges
to $f$ in the uniform metric.

\subsection[A property of the L1-approximant for p less than 1]{A property of the $L^{1}$-approximant for $p<1$}\label{subsec-appl-bl}
We continue the preceding example, but now let $0<p<1$. In 1994, L.G.~Brown
and B.J.~Lucier \cite{BroLuc94} proved the following statement, which can be
regarded as an analogue of Jackson's result \eqref{jac-res}. The approximant
$A_{1}f$ is determined by the criterion for best approximation in $L^{1}$,
which allows the operator $A_{1}$ to be extended to all measurable functions
$f\in L^{0}(M)$ (see \cite{BroLuc94} for details). Then, for $0<p<1$,
\[
\|f-A_{1}f\|_{p}\le (2\mathcal{C}_{1}(Y)+1)^{1/p}\|f-A_{p}f\|_{p}.
\]
A stronger inequality with the constant
$(2\mathcal{C}_{1}(Y))^{1/p-1}$ instead of
$(2\mathcal{C}_{1}(Y)+1)^{1/p}$ was stated in \cite{BroLuc94}. With this
constant the two sides of the result are consistent at $p=1$. However, this
stronger inequality remained unproved.

\subsection{The Remez problem}\label{subsec-appl-remez}
Consider a version of the extremal Remez problem concerning concentration of
the norm of functions from $Y\subset C(M)\cap L^{p}(M)$ on a set of small
measure (see, for example, \cite{MalRyu14,TemTik17,TikYud20}). How small can
the measure of a measurable set $E\subset M$ be if there exists $f\in Y$ such
that
\[
\frac{1}{2}\int_{M}|f|^{p}\,dx\le \int_{E}|f|^{p}\,dx
\]
(the constant $\frac{1}{2}$ can be replaced by another fixed number)? The
following simple argument gives a lower estimate for $|E|$ (see, for example,
\cite{TemTik17}). Continuing the preceding inequality, we obtain
\[
\int_{E}|f|^{p}\,dx\le |E|\,\|f\|_{\infty}^{p}\le
|E|(\mathcal{C}_{p}(Y))^{p}\|f\|_{p}^{p}.
\]
Therefore,
\begin{equation}\label{E-C}
|E|\ge \frac{1}{2(\mathcal{C}_{p}(Y))^{p}}.
\end{equation}

The Remez problem in the uniform metric has been solved for algebraic
polynomials and, only recently, for trigonometric polynomials (see
\cite{TikYud20}). In the $L^{1}$ metric, the problem for trigonometric
polynomials and entire functions of exponential type was solved in
\cite{MalRyu14}. The author is not aware of solutions in other cases. In
particular, \cite{MalRyu14} proves that for $Y=\mathcal{T}_{n}$ in
$L^{1}(\mathbb{T})$, for sufficiently large $n$, the sharp inequality is
$|E|\ge \frac{\pi}{n+1}$ (equality is attained for the Rogosinski polynomial).
By \eqref{E-C} and the results of Section~\ref{sec-hist} (see \eqref{cL} and
$L\approx 0.172$),
\[
|E|\ge \frac{1}{2c_{n}}\ge \frac{1}{2(n+1)L}\approx \frac{2.9}{n+1},
\]
which is, of course, somewhat weaker than the sharp result.

Another example is provided by \cite{TemTik17,TikYud20}. Let
$T\in\mathcal{T}_{n}$ and $s>0$. Choose
\[
E=\{t\in\mathbb{T}\colon |T(t)|>\|T\|_{1}/s\}.
\]
Then $|E|\le s$, and the sharp Remez inequality
$\|T\|_{\infty}\le R(s)\|T\|_{\infty;\mathbb{T}\setminus E}$ implies
\[
\|T\|_{\infty}\le \frac{R(s)}{s}\|T\|_{1}.
\]
As $s\to 0$,
$R(s)=1+\frac{(ns)^{2}}{8}+O(n^{4}s^{4})$. Taking $s=\frac{1}{n}$ and using
the exact formula for the Remez constant, we obtain
\[
\|T\|_{\infty}\le \Bigl(\cosh\frac12+o(1)\Bigr)n\|T\|_{1}.
\]
This has the correct order, but its constant is substantially larger than the
sharp constant $c_{n}\approx 0.172\,n$.

In view of \eqref{E-C}, it is of interest to identify cases in which an upper
estimate for $|E|$ has the same order as the lower estimate given by the
Nikolskii constant. In some cases a good upper estimate is obtained from the
extremal function in Logan's extremal problem
(see~\cite{GorIvaTik20,DaiGorTik21}), where one seeks a function of mean zero
whose positivity set has small measure, as below.

It is also useful to mention related results on estimating the measure of the
positivity set of mean-zero functions from $Y$ (see, for example,
\cite{Yud02,Yud04}). Let $p=1$, $f\in Y\setminus \{0\}$, and
$\int_{M}f\,dx=0$. How small can the measure of
$E=\{x\in M\colon f(x)>0\}$ be? The following argument was used in
\cite{Yud04}:
\[
\int_{M}|f|\,dx=\int_{E}f\,dx-\int_{M\setminus E}f\,dx=2\int_{E}f\,dx\le
2|E|\,\|f\|_{\infty}\le 2|E|\mathcal{C}_{1}(Y)\|f\|_{1},
\]
which gives the analogue of \eqref{E-C} for $p=1$,
\begin{equation}\label{E-C-1}
|E|\ge \frac{1}{2\mathcal{C}_{1}(Y)}.
\end{equation}
Again, one is led to the problem of finding good upper estimates for $|E|$.
This problem is particularly difficult, for example, for the subspace of
spherical harmonics of degree $n$ on the sphere. This subspace is an
eigenspace of the Laplace--Beltrami operator $\Delta_{0}$. The problem is
studied in the broader context of the distribution of nodal domains of
eigenspaces of powers of the Laplacian. Many difficult open problems remain
here (see, for example, \cite{Log18}).

\subsection[L1 approximation]{$L^{1}$ approximation}
In \cite{BenKroPin12}, the authors considered the following question in
$L^{1}$ approximation. In the setting of Subsection~\ref{subsec-appl-remez},
let $f\in L^{1}(M)$ vanish identically outside a set $E\subset M$. How large
can $|E|$ be while the approximant $A_{1}f=0$ is still guaranteed? Let
$\alpha^*(Y)$ denote the largest number for which this property holds for all
such functions whenever $|E|\le\alpha^*(Y)$. Using the criterion for $A_{1}f$
mentioned in Subsection~\ref{subsec-appl-bl}, the authors proved that for
$q>1$,
\[
\alpha^*(Y)\ge \frac{1}{(2\mathcal{C}_{1q}(Y))^{q'}}\ge
\frac{1}{2^{q'}\mathcal{C}_{1}(Y)},
\]
where the second inequality follows from \eqref{pq}. Thus we again obtain an
estimate related to \eqref{E-C} and \eqref{E-C-1}. The difficult problem of a
correct upper estimate also arises here; a partial answer was given in
\cite{BenKroPin12}.

Using the Nikolskii constant $L=\mathcal{C}_{1}(\mathcal{E}_{1})$ as an
example, we show another connection with best weighted $L^{1}$ approximation.
We have
\[
L^{-1}=\inf\,\biggl\{\int_{\mathbb{R}}|f|\,dx\colon f\in
\mathcal{E}_{1}^{1,\text{even}},\ f(0)=1\biggr\}.
\]
Represent $f$ in the form $f(x)=1-x^{2}g(x)$, where
$g\in\mathcal{E}_{1}^{\text{even}}$ and
$x^{-2}-g\in L_{x^{2}}^{1}(\mathbb{R})$. Then
\[
\int_{\mathbb{R}}|f(x)|\,dx=\int_{\mathbb{R}}|1-x^{2}g(x)|\,dx=
\int_{\mathbb{R}}|x^{-2}-g(x)|\,x^{2}\,dx.
\]

Thus, computing $L$ reduces to the problem of best $L^{1}$ approximation of
$x^{-2}$ by entire functions of exponential type at most $1$ in the weighted
space $L_{x^{2}}^{1}(\mathbb{R})$ (see \cite{VinGla14,LitSpa18}). The
function $x^{-2}$ itself does not belong to this space, which is not an
obstacle \cite{VinGla14}. In \cite{VinGla14,LitSpa18}, the problem of
approximating individual functions by entire functions of exponential type at
most $\sigma$ was solved in $L_{|E|^{-2}}^{1}(\mathbb{R})$, where $E$ is a
Hermite--Biehler entire function of exponential type $\sigma/2$ \cite{Lev56}.
However, as noted in \cite{LitSpa18}, the weight $x^{2}$ does not fit into this
scheme, so this approach does not yet yield a solution of the Nikolskii
problem for~$L$.

\subsection{Number theory}\label{subsec-appl-num}
In the recent paper \cite{CarMilSou19}, E.~Carneiro, M.~B.~Milinovich, and
K.~Soundararajan proved that, assuming the Riemann hypothesis for the zeros of
the zeta function, the prime numbers $p_n$ satisfy
\begin{equation}\label{p-p}
\limsup_{n\to \infty}\frac{p_{n+1}-p_{n}}{\sqrt{p_{n}}\log p_{n}}\le \frac{21}{25}.
\end{equation}
The existence of a finite constant on the right-hand side had been proved by
H.~Cram\'er in 1920, and in 2015 A.~Dudek showed that the constant is at most
$1$. Several related number-theoretic results were also proved in
\cite{CarMilSou19}. The authors used the following one-dimensional Fourier
optimization problems, one of which reduces to the constant $L$ (see above
and Subsection~\ref{subsec-hist-now}). Let the Fourier transform be defined by
$\widehat{F}(t)=\int_{-\infty}^{\infty}e^{-2\pi ixt}F(x)\,dx$. For
$1\le A<\infty$, one has to determine
\begin{align*}
\mathcal{C}(A)&=\sup_{F\in C(\mathbb{R})\cap L^{1}(\mathbb{R})}
\frac{1}{\|F\|_{1}} \biggl(|F(0)|-A\int_{\mathbb{R}\setminus
[-1,1]}|\widehat{F}(t)|\,dt\biggr),\\
\mathcal{C}^{+}(A)&=\sup_{F\in C^{\text{real}}(\mathbb{R})\cap
L^{1}(\mathbb{R})} \frac{1}{\|F\|_{1}} \biggl(F(0)-A\int_{\mathbb{R}\setminus
[-1,1]}(\widehat{F}(t))_{+}\,dt\biggr).
\end{align*}
The constant $\frac{21}{25}$ in \eqref{p-p} follows from the estimate
$\mathcal{C}^{+}(\frac{36}{11})>\frac{25}{21}$. In the limit
$A\to\infty$, one obtains the Nikolskii constant
$\mathcal{C}(\infty)=2\pi L$.

\subsection{Metric geometry}\label{subsec-appl-geom}
In Subsection~\ref{subsec-common-shift} we considered the shifted Nikolskii
problem for nonnegative functions
\[
\mathcal{A}_{1}^{+}(Y_{n0};0)=\sup_{0\le f\in Y_{n0}\cap
L^{1}(M)}\frac{f(0)}{\|f\|_{1}}
\]
and showed how it can be solved by quadrature formulas. This problem and its
variants are also actively studied in harmonic analysis because of their
applications, for example, in metric geometry.

Let us give some examples. Let $M=\mathbb{S}^{d}$. Then the results of
Subsection~\ref{subsec-common-pos} (see also \cite{DaiGorTik20-anal}) give
\[
|\mathbb{S}^{d}|\mathcal{A}_{1;\mathbb{S}^{d}}^{+}(\Pi_{n}^{d};e_{d+1})=
\mathcal{A}_{1,w_{d};[-1,1]}^{+}(\mathcal{P}_{n}^{\text{real}};1)=\sup_{0\le
P\in \mathcal{P}_{n}}\frac{P(1)}{\int_{-1}^{1}P(t)w_{d}(t)\,dt}=I_{d,n},
\]
where $e_{d+1}$ is the north pole of the sphere and
$w_{d}=\frac{(1-t^{2})^{d/2-1}}{\int_{-1}^{1}(1-t^{2})^{d/2-1}\,dt}$ is the
normalized algebraic Gegenbauer weight. It turns out that the solution of the
extremal problem $I_{d,n}$ for polynomials gives the well-known
Delsarte--Goethals--Seidel bound for tight designs \cite{DelGoeSei77}. A finite
set of points $x_{1},\ldots,x_{N}\in \mathbb{S}^{d}$ is called an $n$-design
if the quadrature formula
\[
\frac{1}{|\mathbb{S}^{d}|}\int_{\mathbb{S}^{d}}f(x)\,dx=\frac{1}{N}\sum_{i=1}^{N}f(x_{i})
\]
holds for every algebraic polynomial of degree at most $n$ in $d+1$ variables
(as in \eqref{poly-d}). A fundamental problem in metric geometry (and in its
applications, for example, to coding theory) is to find $n$-designs with the
smallest possible number of nodes $N(d,n)$. It is known that
\cite{DelGoeSei77,Lev98}
\[
N(d,n)\ge I_{d,n}= \binom{d+[\frac{n+1}{2}]-1}{d}+\binom{d+[\frac{n}{2}]}{d},
\]
which is the tight-design bound. This estimate was extended to compact
rank-one Riemannian manifolds in \cite{Lev98}.

Let $M=\mathbb{R}^{d}$ and let $\Omega\subset \mathbb{R}^{d}$ be, as above, a
convex body. Then
\[
\mathcal{A}_{1;\mathbb{R}^{d}}^{+}(\mathcal{E}_{1}(\Omega);0)= \sup_{0\le f\in
\mathcal{E}_{1}(\Omega)\cap
L^{1}(\mathbb{R}^{d})}\frac{f(0)}{\int_{\mathbb{R}^{d}}f(x)\,dx}=I(\Omega),
\]
and, after passing to the Fourier transform,
\begin{equation}\label{I-prob}
I(\Omega)=\sup\,\Bigl\{\frac{\widehat{f}(0)}{f(0)}\colon f\in
C(\mathbb{R}^{d})\cap L^{1}(\mathbb{R}^{d}),\ \operatorname{supp}f\subset
\Omega,\ \widehat{f}\ge 0\Bigr\}.
\end{equation}

This is the so-called Tur\'an extremal problem for positive definite functions
with small support. Many papers are devoted to this problem and its variants,
including joint papers by the author (see, for example,
\cite{GorTik18,GorTik19,GorIva19-sb} and the references therein). At present,
the problem $I(\Omega)$ has been solved only for space-tiling bodies (Voronoi
polytopes and spectral bodies) and for the Euclidean ball, with proofs by
several authors \cite{Sie35,Vaa85,Gor01,AreBer02,KolRev03,BiaKel15}. In these
cases
$I(\Omega)=\frac{|\frac{1}{2}\,\Omega|}{(2\pi)^{d}}$, and the extremal
function is the convolution
$\chi_{\frac{1}{2}\,\Omega}*\chi_{\frac{1}{2}\,\Omega}$, whose Fourier
transform corresponds to the Fej\'er kernel. A difficult open problem is to
prove that this function is extremal for other convex bodies $\Omega$.

If in \eqref{I-prob} we set $\Omega=B^{d}$ and enlarge the class of admissible
functions by allowing them to take nonpositive values outside $B^{d}$, we
obtain the well-known linear-programming bound for the density of sphere
packings in $\mathbb{R}^{d}$ by equal balls of radius $1/2$ (see
\cite{Gor00,CohElk03}). The corresponding problem is also called the Delsarte
problem. Recently, it was solved by means of modular forms and Schwartz
functions for $d=8,24$ \cite{Via17,CohKumMilRadVia17}, yielding in these
dimensions a solution of the higher-dimensional Kepler problem for sphere
packing. This outstanding result emphasizes the importance of studying
variants of Nikolskii-type problems. The theory of the Delsarte problem is
developing in several directions. One of them is related to the uncertainty
principle of J.~Bourgain, L.~Clozel, and J.-P.~Kahane \cite{BouCloKah10}, where
recent progress was also obtained in \cite{CohGon18,GonOliRam21}, as well as
in \cite{GorIvaTik20} for functions with bounded spectrum.

\section{Conclusion}

We conclude this survey with several open problems in the theory of sharp
Bernstein--Nikolskii inequalities that arose in the discussion above. Some of
them are difficult, while others are more specialized, but each is important
in its own way. We do not repeat the definitions; see the preceding sections.

\begin{enumerate}
\item Prove a Levin--Lubinsky type result for $q\ne \infty$.

\item Give examples of the maximum-shift method in the weighted case for
$p<1$.

\item In \eqref{C-P-cheb}, determine the asymptotic behavior of
$\|P_{n}^{*}\|_{p,w_{\alpha,\beta}}$ as $n\to \infty$.

\item Prove estimates of type \eqref{n-C-l-n} for
$\mathcal{T}_{n}(B^{d})$.

\item Compute the constant $L=\mathcal{L}_{11}$, for example, as a root of an
explicit transcendental equation.

\item Can the upper and lower bounds in \eqref{LA-d} for the multidimensional
constant $\mathcal{L}_{d1}$ be improved exponentially in $d$?

\item Prove that the radial extremal function in the problem
$\mathcal{L}_{d1}$ is unique not only among radial functions.

\item Prove that the signs of the nonzero Taylor coefficients of extremal
functions in the problems $c_{np}^{(r)}$ and $L_{p}^{(r)}$ alternate.

\item Solve problem \eqref{j-j} for $-1/2<\alpha<-0.272$.

\item Develop the theory presented in Section~\ref{sec-common} for hyperbolic
spaces.

\end{enumerate}

\end{document}